\documentclass[preprint,11pt]{elsarticle}

\usepackage{amssymb}
\usepackage{amsmath}
\usepackage{hyperref}

\usepackage{physics}
\usepackage{tensor}
\usepackage{rotating}
\usepackage[toc,page]{appendix}
\usepackage{csquotes}
\usepackage{dsfont}
\usepackage{mathrsfs}
\usepackage{centernot}
\usepackage{pdfpages}
\usepackage{float}
\usepackage{orcidlink}

\usepackage[margin=1in]{geometry}

\journal{Differential Geometry and Its Applications}

\begin{document}

\begin{frontmatter}


% \title{Rényi-Induced Information Geometry and Hartigan's Prior Family\tnoteref{label1}}

\author{Rebecca Maria Kuntz\fnref{label1}\orcidlink{0009-0006-0960-817X}}
\ead{kuntz@stud.uni-heidelberg.de}
\author{Heinrich von Campe\fnref{label2,label3}\orcidlink{0009-0006-4071-3576}}
\ead{heinrich.campe@iwr.uni-heidelberg.de}
\author{Björn Malte Schäfer\fnref{label1}\orcidlink{0000-0002-9453-5772}}
\ead{bjoern.malte.schaefer@uni-heidelberg.de}


\title{Rényi-Induced Information Geometry and Hartigan's Prior Family}

% \author[label1]{Rebecca Maria Kuntz}
% \author[label2]{Heinrich von Campe}
% \author[label1]{Björn Malte Schäfer}

\affiliation[label1]{organization={Zentrum für Astronomie der Universität Heidelberg, Astronomisches Rechen-Institut},%Department and Organization
            addressline={\newline Philosophenweg~12}, 
            city={Heidelberg},
            postcode={69120}, 
            country={Germany}}

\affiliation[label2]{organization={Interdisziplinäres Zentrum für wissenschaftliches Rechnen der Universität Heidelberg},
            addressline={\newline Im Neuenheimer Feld~205},
            city={Heidelberg},
            postcode={69120},
            country={Germany}}

\affiliation[label3]{organization={Zuse School ELIZA}}

\newcommand{\orcidauthorA}{0009-0006-0960-817X}% Add \orcidA{} behind the author's name
\newcommand{\orcidauthorB}{0009-0006-4071-3576}
\newcommand{\orcidauthorC}{0000-0002-9453-5772}

% \author{Rebecca Maria Kuntz, Heinrich von Campe and Björn Malte Schäfer}%% Author name

% \affiliation{organization={ Affiliation: Zentrum für Astronomie der Universität Heidelberg, Astronomisches Rechen-Institut},%Department and Organization
%             addressline={Philosophenweg 12}, 
%             city={Heidelberg},
%             postcode={69120}, 
%             country={Germany}}

%% Abstract
\begin{abstract}
We derive the information geometry induced by the statistical Rényi divergence,
namely its \mbox{metric} tensor, its dual parametrized connections, as well as its dual
Laplacians. Based on these results, we demonstrate that the Rényi-geometry, though closely related, differs in structure from Amari's well-known $\alpha$-geometry. Subsequently, we derive the canonical uniform prior distributions for a statistical manifold endowed with a Rényi-geometry, namely the dual Rényi-covolumes. We find that the Rényi-priors can be made to coincide with
Takeuchi and Amari's $\alpha$-priors by a reparameterization, which is itself of particular significance in statistics. Herewith, we demonstrate that Hartigan's parametrized ($\alpha_H$) family of priors is precisely the parametrized ($\rho$) family of Rényi-priors $(\alpha_H = \rho)$.
\end{abstract}


% \begin{highlights}
% \end{highlights}

%% Keywords

\end{frontmatter}

\section{Introduction}\label{sect: introduction}
The field of information geometry studies the innate geometry of statistics, opening a new, intuitive way to reason about the geometric invariants of statistical model families. In contrast to the geometries induced by Amari's $\alpha$-divergences \cite{amari_information_2016}, the geometry induced by the statistical Rényi divergence, whose associated entropy fulfills the generalized Shannon-Khinchin axioms \cite{Shannonkhinchin0, Shannonkhinchin}, remains largely unstudied. Motivated by this gap, this work constructs the Rényi-geometry, key quantities of which are contextualized with statistical estimators so as to discern their significance
for statistical inference. Evaluation of the Rényi-geometry against the backdrop of Amari's \mbox{$\alpha$-geometry \cite{amari_information_2016}} shows their inherent difference in structure. Dual affine volume forms are derived for the Rényi-geometry, leading to the novel Rényi-priors, which are shown to precisely coincide with Hartigan's parametrized prior family. This result offers a new geometric reasoning for a result by Takeuchi and Amari \cite{alphaparallelpriors}.

Previous works on geometries related to the Rényi divergence include de Souza, Vigelis and Cavalcante \cite{Souza2016}, who use the Rényi divergence to study a generalization of the $\alpha$-geometry. Furthermore, van Erven and Harremo{\"e}s \cite{van_erven_renyi_2014} give a review of the Rényi divergence and its properties. Studies addressing divergence-induced canonical prior choices in information geometry include those by Takeuchi and Amari \cite{alphaparallelpriors} and Jiang, Tavakoli and Zhao \cite{weylprior}.

The functional idea of information geometry begins with a sample (data) space $\mathscr{Y}$ endowed with a suitable integration
measure $\dd \mu$. Furthermore, let $p$ denote a probability
distribution (statistical model),
\begin{equation}
  p : \Theta \to \mathcal{P}(\mathscr{Y}).
\end{equation}
on $\mathscr{Y}$. Here, $\Theta \subseteq \mathds{R}^n$ is referred to as the parameter space. $\mathcal{P}(\mathscr{Y})$ denotes a certain set of such models. Subsequently, an information
geometry can be constructed on the topological space constituted by such a model
family $\mathcal{P}$.

The most prominent examples of such families are the
exponential family \cite{amari_information_2016}
\begin{equation}
    \mathcal{P}_e = \Big\{ p_\theta(y) = \exp\left( \theta^i y_i + k(y)-\psi(\theta)\right)\Big\}_{\theta \in \Theta}\,,\label{eq: exponential family}
 \end{equation}
with $k : \mathscr{Y} \rightarrow \mathds{R} $ a function of the data $y$, $\psi: \Theta \rightarrow  \mathds{R} $ a normalization function,
as well as the mixture family \cite{amari_information_2016,Nielsen_2020}
  \begin{equation}
    \mathcal{P}_m = \Big\{p_\eta(y) = \eta^i F_i(y)+ C(y) \Big\}_{\eta\in \Theta}\,,\label{eq: mixture family}
  \end{equation} 
with $F_i: \mathscr{Y} \to \mathds{R}$ some linearly independent functions such that $\smallint_{\mathscr{Y}} F_i(y) = 0$ and $C: \mathscr{Y} \to \mathds{R}$ a function s.t.
$\smallint_{\mathscr{Y}} C(y) = 1$.

The key insight of information geometry is that families of probability
distributions such as the above may be described as so-called statistical
manifolds \cite{amari_information_2016,Nielsen_2020}. In this picture, the parameters
locally provide a coordinate chart \cite{giesel_information_2021}. This is a generalization of the more classical picture
where model parameters live on $\mathds{R}^n$. An illustration
of this interpretation is that the reparametrization of a model $p$ is now
interpreted as a mere change of coordinates on the corresponding statistical
manifold $\mathcal{M}$ \cite{Nielsen_2020}.

The aim of information geometry is to define
geometric structures on these manifolds in a way that illuminates the inherent
geometric properties of the statistical models. Rao's highly innovative work on the Kullback-Leibler (KL) divergence \cite{rao1945information} uncovered the key role played by the the Fisher information metric \cite{original_fisher_information} in information geometry,
\begin{equation}
\tensor[^{\text{F}}]{g}{_i_j}(\theta) := \int \dd{y} p \, \partial_i \log p \, \partial_j \log p \,,\label{eq: Fisher definition}
\end{equation}
since it promotes the statistical manifold to a metric space and is invariant under sufficient statistics~\cite{amari_information_2016}.

The Kullback-Leibler divergence in conjunction with the Fisher information is widely used in machine learning and
statistical inference for methods such as Natural Gradient Descent \cite{natural_gradient_descent} 
and Riemann Manifold
Hamilton Monte Carlo \cite{rmhmc}.

Given this metric tensor, Efron and Dawid
\cite{efron} found that the most natural construction
from an information geometric perspective are \emph{two} connections (the
so-called $e$ and $m$ connections) which are \emph{dually} metric compatible. It was Amari who recognized that information geometry is thus dual in essence \cite{amari_information_2016}. Building on
this, canonical dual prior distributions can be derived from volume forms which
are \emph{parallel} to (preserved under) such dual connections.

\begin{figure}[H]
  \centering 
  \def\svgwidth{0.8\textwidth}
  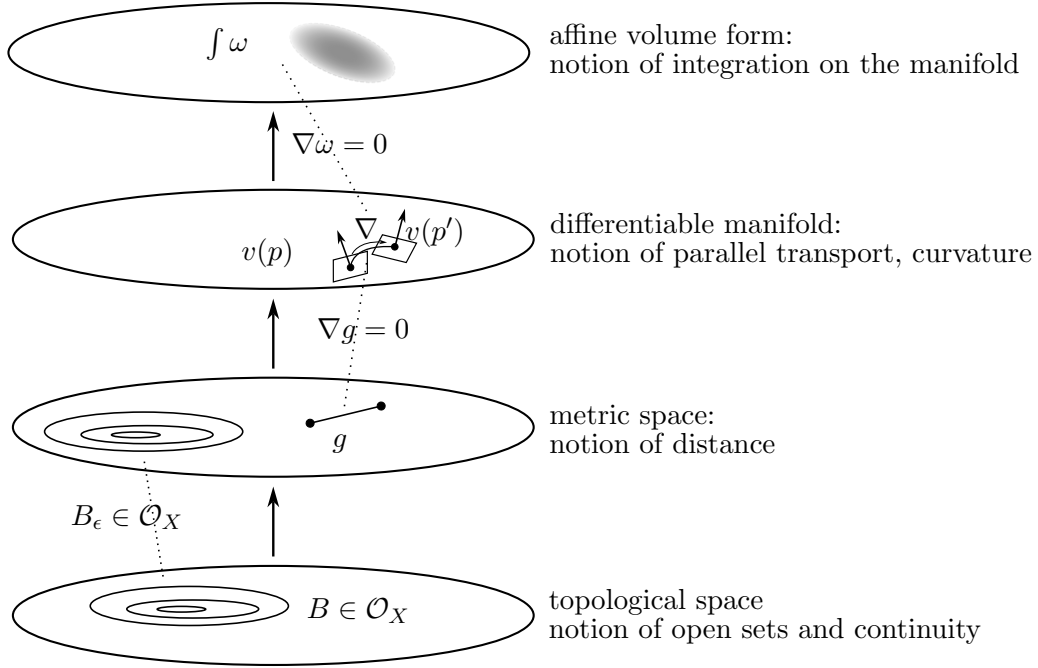
  \caption{Geometric structures are succesively built on a topological space.}
  \label{fig: cake_layers}
\end{figure}

The most famous example of such geometric priors is Jeffreys prior,
which follows from the volume form parallel to the Levi-Civita connection with
the Fisher metric \cite{jeffreysprior101}. \autoref{fig: cake_layers} schematically depicts the construction of different layers of geometrical structure on a topological space.
A natural way to construct the aforementioned geometry on a statistical manifold
is from statistical divergences. Loosely speaking, a statistical divergence or relative entropy is a
macroscopic dissimilarity measure between probability distributions \cite{Nielsen_2020,manyfaces}. A canonical choice for such a divergence is the KL divergence, whose associated entropy, the Shannon entropy, is axiomatically singled out by the so-called Shannon-Khinchin axioms \cite{Shannonkhinchin0}. The generalized Shannon-Khinchin axioms lead to the Rényi entropy as a generalization of the Shannon entropy, making the Rényi divergence a suitable choice for inducing an information geometry \cite{Shannonkhinchin}.

Eguchi \cite{eguchi1983divergenceinduced} demonstrates that any such divergence permits the definition of a
metric and a pair of dual connections. The overall geometric structure of the manifold thus not only depends on
the statistical model at hand but also the statistical divergence that induces different geometric objects.

This paper is organized as follows: Section \ref{sect: Rényi Metric and
Connections} discusses several exisiting divergence-induced geometries before presenting the Rényi divergence's induced geometry.  Section~\ref{sect: Dual Rényi-Laplace-Beltrami Operators} derives the dual
Laplace-Beltrami operators for the Rényi-geometry. Subsequently, \mbox{Section
\ref{sect: Rényi Volume Forms, Rényi-Priors}} reviews various existing canonical
priors of information geometry before deriving the Rényi-priors from volume forms
parallel to the dual Rényi-connections which are shown to coincide with Hartigan's prior family.
Finally, Section \ref{sect: Conclusions and Outlook} summarizes our results and
gives a synopsis of possible next steps.

\section{Rényi Metric and Connections}\label{sect: Rényi Metric and Connections}
\subsection{Preliminaries I: (Dual) Connections on (Statistical) Manifolds} \label{sec: dual connections}

Given a manifold $\mathcal{M}$, a connection $\nabla : \Gamma(T\mathcal{M}) \times \Gamma(T\mathcal{M}) \to
\Gamma(T\mathcal{M})$
is a covariant generalization of the directional derivative. Here,  $\Gamma(T\mathcal{M})$ is the associated tangent bundle. The connection is bilinear and fulfills
\cite{Bartelmann_2019}
\begin{equation}
\nabla_{f V} Y = f \nabla_V Y\,, \quad
\nabla_V (f Y) = f \nabla_V Y + V(f) Y\,,
\end{equation}                
for any vector fields $V,Y \in \Gamma(T\mathcal{M})$ and any smooth function 
$f \in C^{\infty}(\mathcal{M})$. It may be fully characterized 
by its \emph{connection coefficients}
$(\nabla_{\partial_i} (\partial_j))^k
= - (\nabla_{\partial_i} \dd{x}^k)_j
= \tensor[]{\Gamma}{^k_{ij}}$.

If the manifold $\mathcal{M}$ is equipped with a metric $g$, we may furthermore define the
notion of \emph{metric compatibility} of an affine connection $\nabla$
w.\,r.\,t.\@ the metric $g$ by
\cite{Bartelmann_2019}
\begin{equation} \label{eq: metric compatibility}
Z g(X,Y)
= g(\tensor{\nabla}{_Z} X, Y) + g(X, \tensor{\nabla}{_Z} Y)\,,
\end{equation}
for arbitrary  vector fields $X,Y,Z \in \Gamma(T\mathcal{M})$. This
equation directly implies that the the inner product of two vectors $g(X,Y)$ that are
parallel transported along a curve $\gamma$ by means of a connection $\nabla$ (i.\,e.\@ $\nabla_{\dot{\gamma}} X = 0$) will remain
constant if the connection is metric compatible.

One commonly uses the Levi-Civita connection
$\tensor[^{\text{LC}}]{\nabla}{}$
which is the unique torsion-free metric-compatible connection
\cite{Bartelmann_2019}. Information geometry on the other hand is built on the notion of \emph{dual metric
compatibility} of \emph{two} connections $\nabla, \nabla^*$ w.\,r.\,t.\@ a
metric $g$, defined by
\cite{Nielsen_2020}
\begin{equation} \label{eq: dual metric compatibility}
Z g(X,Y)
= g(\nabla_Z X, Y) + g(X, \nabla^*_Z Y)\,.
\end{equation}
In general, neither $\nabla$ nor $\nabla^*$ are metric compatible by
themselves.  In analogy to the description above, \autoref{eq: dual metric compatibility} indicates that the inner
product of two vectors $g(X,Y)$ remains constant along a curve $\gamma$ if
one of the vectors is parallel transported by means of $\nabla$ and the other
with $\nabla^*$  \cite{Nielsen_2020}. One may easily see that the arithmetic mean of the
two is indeed the Levi-Civita connection,
$\tensor[^{\text{LC}}]{\nabla}{} = \frac{1}{2} (\nabla + \nabla^*)$.
The difference between different connection coefficients is always a
tensor~\cite{schaefer_tooltips_2022}.
Connections allow to study the notion of curvature of 
manifolds~\cite{Bartelmann_2019}. Note that a manifold is called flat (has zero curvature) if there exists a coordinate choice for which the connection coefficients vanish globally.

\subsection{Preliminaries II: Divergence-Induced Information Geometries}\label{sect: Divergence-induced geometries on statistical manifolds}
Given a statistical manifold $\mathcal{M}$ and a statistical divergence 
\mbox{$D: \mathcal{M} \times \mathcal{M} \to \mathds{R}_{+} \cup \{0\}$} with
\begin{equation}
D[\theta : \theta'] = 0 \leftrightarrow \theta = \theta'\,, \quad
\partial_i D[\theta : \theta'] |_{\theta = \theta'}
= 0 =
\partial_{j'} D[\theta : \theta'] |_{\theta = \theta'} \, ,
\end{equation}
\begin{equation}
-\partial_i \partial_{j'} D[\theta : \theta'] |_{\theta = \theta'} \quad
\text{is positive definite}\,,
\end{equation}
Eguchi \cite{eguchi1983divergenceinduced} demonstrated that one  may construct a
divergence-induced information geometry on~$\mathcal{M}$ with the following quantities:
\begin{align}
\tensor[^D]{g}{_{ij}} = -\partial_i \partial_{j'} D[\theta : \theta']\big|_{\theta = \theta'}\,, \quad
&\text{metric tensor}, \label{eq: metric from divergence} \\
\tensor[^D]{\Gamma}{_{ijk}} = -\partial_i \partial_j \partial_{k'} D[\theta : \theta']\big|_{\theta = \theta'}\,, \quad
&\text{connection coefficient of first kind}, \label{eq: connection from divergence} \\
\tensor[^{D^*}]{\Gamma}{_{ijk}} = -\partial_k \partial_{i'} \partial_{j'} D[\theta : \theta']\big|_{\theta = \theta'}\, , \quad
&\text{dual connection coefficient of first kind}. \label{eq: dual connection from divergence} 
\end{align}
One may easily check  that the behaviour of these quantities under coordinate
transformations on $\mathcal{M}$ \cite{Bartelmann_2019}
is correct. 
Please note that
these are the connection coefficients of the first kind which are related to the
aforementioned connection coefficients of the second kind by
$\Gamma^k_{ij} = g^{k\ell} \Gamma_{ij\ell}$ \cite{Bartelmann_2019}.
Furthermore, note that the connections  
$ \tensor[^D]{\nabla}{}, \tensor[^{D^*}]{\nabla}{} $
induced from a divergence $D$ will be always torsion-free,
$ \tensor[^{D^{(*)}}]{\Gamma}{_{ijk}} = \tensor[^{D^{(*)}}]{\Gamma}{_{jik}}$,
since the derivatives in 
\autoref{eq: connection from divergence} and \autoref{eq: dual connection from divergence}
commute.

It is straightforward to
see that divergence-induced connections automatically fulfill
dual metric compatibility, 
\begin{equation}
\partial_i \tensor[^D]{g}{_j_k}
= -\partial_i \partial_j \partial_{k'} D[\theta : \theta']\big|_{\theta = \theta'}
-\partial_{i'} \partial_j \partial_{k'} D[\theta : \theta']\big|_{\theta = \theta'}
=  \tensor[^D]{\Gamma}{_{ijk}} + \tensor[^{D^*}]{\Gamma}{_{ikj}}\,,
\end{equation}
which implies \autoref{eq: dual metric compatibility} in local coordinates.
As for all connections, the difference between 
$ \tensor[^D]{\nabla}{}$ and $ \tensor[^{D^*}]{\nabla}{} $
will be always tensorial, at the same time their average will indeed be the
Levi-Civita connection w.\,r.\,t.\@ the metric tensor $\tensor[^D]{g}{}$.

\subsection{Kullback-Leibler-Geometry} \label{sect: KL geometry}
One instance of such a geometric structure is the geometry induced by the
Kullback-Leibler divergence~\cite{kullback_leibler_original_definition},
\begin{equation}
D_{\text{KL}}[\theta : \theta'] = \int \dd \mu(y) \, p \log \frac{p}{p'}.
\end{equation}
Using \autoref{eq: metric from divergence} -- \autoref{eq: dual connection from divergence},
we find
\begin{align}
\tensor[^{\text{KL}}]{g}{_{ij}}
&= \int \dd \mu(y) p \, \partial_i \log p \, \partial_j \log p\,
= \tensor[^{\text{F}}]{g}{_{ij}}\,, \\
\tensor[^{\text{KL}}]{\Gamma}{_{ijk}}
&= \int \dd \mu(y) p \,\partial_i \partial_j \log p \,\partial_k \log p\,
+ \int \dd \mu(y) p \,\partial_i \log p \,\partial_j \log p \,\partial_k \log p \\
&=: \tensor[^{(e)}]{\Gamma}{_{ijk}} + C_{ijk}
=: \tensor[^{(m)}]{\Gamma}{_{ijk}}
\label{eq: mixture connection coefficient and Amari Chentsov defintion}
\hspace{6.3em}  \leftrightarrow \quad 
\tensor[^{\text{KL}}]{\nabla}{} = \tensor[^{(m)}]{\nabla}{} \, , \\
\tensor[^{\text{KL}^*}]{\Gamma}{_{ijk}}
&= \int \dd \mu(y) p \,\partial_i \partial_j \log p \,\partial_k \log p\,
=: \tensor[^{(e)}]{\Gamma}{_{ijk}}
\label{eq: exponential connection coefficient defintion}
\quad  \leftrightarrow \quad
\tensor[^{\text{KL}^*}]{\nabla}{} = \tensor[^{(e)}]{\nabla}{} \, .
\end{align}
The metric may be identified as the Fisher metric. It is a well fact in
information geometry~\cite{Nielsen_2020}
and uniquely characterized by its invariance under
sufficient statistics, that is, transformations of the data space that preserve the
data's information content about the parameters~\cite{amari_information_2016,Nielsen_2020}.
The connections 
$\tensor[^{\text{KL}}]{\nabla}{}$
induced by the KL divergence and its dual
$\tensor[^{\text{KL}^*}]{\nabla}{}$ may be identified
as the connections of the mixture and exponential family,
$\tensor[^{(m)}]{\nabla}{}$ and $\tensor[^{(e)}]{\nabla}{}$,
respectively. They bear these names since the corresponding geometry is flat
for the exponential and mixture family $\mathcal{P}_e$ and $\mathcal{P}_m$,
respectively 
\cite{amari_information_2016}
as one may easily see by substituting their definitions 
\autoref{eq: exponential family} and \autoref{eq: mixture family}
in the above expressions for the coefficients. The two connection coefficients
differ by the Amari-Chentsov tensor~\cite{chentsov1982statiscal}
\begin{equation} \label{eq: Amari Chentsov Tensor definition}
C = \tensor[^{(m)}]{\nabla}{} - \tensor[^{(e)}]{\nabla}{},
\end{equation}
which is again uniquely characterised
by its invariance under sufficient statistics
\cite{amari_information_2016}.
Since they are divergence-induced, 
$\tensor[^{(m)}]{\nabla}{}$
and
$\tensor[^{(e)}]{\nabla}{}$
are dually metric compatible w.\,r.\,t.
$\tensor[^{\text{F}}]{g}{}$
\cite{amari_information_2016} and their average yields the  corresponding
Levi-Civita connection,
\begin{equation} \label{eq: Levi-Civita connection definition}
\tensor[^{\text{LC}}]{\nabla}{} 
= \frac{1}{2}  \left( \tensor[^{(m)}]{\nabla}{} + \tensor[^{(e)}]{\nabla}{} \right)
= \tensor[^{(m)}]{\nabla}{} - \frac{1}{2} C
= \tensor[^{(e)}]{\nabla}{} + \frac{1}{2} C \, .
\end{equation}

\subsection{Amari's $\alpha$-Geometry} \label{sect: Amari's alpha geometry}
The first, well-studied generalisation of this construction was done by Amari
\cite{amari_alpha_geometry_original_definition,amari_information_2016}
and Chentsov~\cite{chentsov1982statiscal}
and may be derived by considering the divergence
\begin{equation}
D_\alpha[\theta : \theta'] = \frac{4}{1-\alpha^2} \left(
1 - \int \dd \mu(y) p^\frac{1-\alpha}{2} (p')^\frac{1 + \alpha}{2}
\right),
\end{equation}
where $\rho \in \mathds{R} \setminus \{-1,1\}$. In this case, by employing 
\autoref{eq: metric from divergence} -- \autoref{eq: dual connection from divergence}
we find
\begin{align}
\tensor[^{(\alpha)}]{g}{_{ij}} &= \tensor[^{\text{F}}]{g}{_{ij}}, \\
\tensor[^{(\alpha)}]{\Gamma}{_{ijk}} &= \tensor[^{(e)}]{\Gamma}{_{ijk}} + \frac{1-\alpha}{2} C_{ijk}
\quad \leftrightarrow\quad \hspace{.4em} 
\tensor[^{(\alpha)}]{\nabla}{} = \tensor[^{(e)}]{\nabla}{} + \frac{1-\alpha}{2} C \, ,\label{eq: scaling of alpha connection with C} \\
\tensor[^{(\alpha^*)}]{\Gamma}{_{ijk}} &= \tensor[^{(e)}]{\Gamma}{_{ijk}} + \frac{1+\alpha}{2} C_{ijk}
\quad \leftrightarrow\quad
\tensor[^{(\alpha^*)}]{\nabla}{} = \tensor[^{(e)}]{\nabla}{} + \frac{1+\alpha}{2} C\,.
\end{align}
Firstly, note that by taking the limit of $\alpha \to -1$, we recover the Kullback-Leibler case:
\begin{equation}
\lim_{\alpha \to -1} D_\alpha
= D_{\text{KL}} \, ,
\quad 
\lim_{\alpha \to -1} \tensor[^{(\alpha^{*})}]{\nabla}{}
= \tensor[^{\text{KL}}]{\nabla}{}
\quad \text{and} \quad 
\lim_{\alpha \to -1} \tensor[^{(\alpha^*)}]{\nabla}{}
= \tensor[^{\text{KL}^{*}}]{\nabla}{}.
\end{equation}
While the parameter $\alpha$ does not affect the metric,  it introduces a
continuous tuning mechanism between  the connection 
$\tensor[^{(\alpha)}]{\nabla}{}$ and its dual $\tensor[^{(\alpha^*)}]{\nabla}{}$ 
by scaling the Amari-Chentsov tensor~$C$. Indeed,
$\tensor[^{(-\alpha)}]{\nabla}{} = \tensor[^{(\alpha^*)}]{\nabla}{}$
and
$\tensor[^{(\alpha = 0)}]{\nabla}{} = \tensor[^{\text{LC}}]{\nabla}{}$. Again it is straightforward to see
$\tensor[^{(\alpha)}]{\nabla}{}$ and $\tensor[^{(\alpha^*)}]{\nabla}{}$
fulfill  dual metric compatibility
\cite{amari_information_2016}.

\subsection{Rényi's $\rho$-Geometry}\label{sect: Rényi's rho geometry}
While the $\alpha$-geometry has been well studied in the past, the
geometry induced by the Rényi divergence\footnote{Note: To avoid confusion with
Amari's divergence, we will denote Rényi’s divergence and all derived quantities
with $\rho$ (as in \textbf{R}ényi) instead of $\alpha$.}
\cite{renyi_original,van_erven_renyi_2014},
\begin{equation}
D_\rho[\theta : \theta'] = \frac{1}{\rho - 1} \log \int \dd \mu(y) p^\rho (p')^{1-\rho}\,, \label{eq: renyi divergence}
\end{equation}
with $\rho \in \mathds{R}_+ \setminus \{1\}$, has not yet attracted much attention, despite the fact that the Rényi entropy fulfills the generalized Shannon-Khinchin axioms \cite{Shannonkhinchin}. Now,
\autoref{eq: metric from divergence} -- \autoref{eq: dual connection from divergence}
lead to 
\begin{align}
\tensor[^{(\rho)}]{g}{_{ij}} &= \rho \cdot \tensor[^{\text{F}}]{g}{_{ij}}, \label{eq: renyi metric}\\
\tensor[^{(\rho)}]{\Gamma}{_{ijk}} 
&= \rho \cdot \tensor[^{(e)}]{\Gamma}{_{ijk}} + \rho^2 \cdot C_{ijk}
\quad \quad \quad \hspace{.7em} \leftrightarrow\quad \hspace{.3em}
\tensor[^{(\rho)}]{\nabla}{} = \tensor[^{(e)}]{\nabla}{} + \rho\, C\, , \label{eq: rho connection in terms of e connection and C}\\
\tensor[^{(\rho^*)}]{\Gamma}{_{ijk}} 
&= \rho \cdot \tensor[^{(e)}]{\Gamma}{_{ijk}} + \rho(1-\rho) \cdot C_{ijk}
\quad \leftrightarrow\quad
\tensor[^{(\rho^*)}]{\nabla}{} = \tensor[^{(e)}]{\nabla}{} + (1-\rho)\, C\,.\label{eq: rho* connection in terms of e connection and C}
\end{align}
Note that for the coordinate-free expressions on the right in \autoref{eq: rho
connection in terms of e connection and C} and \autoref{eq: rho* connection in
terms of e connection and C} the connection coefficients of the second kind were
used ($\nabla_i X^j := \partial_i X^j + \tensor{\Gamma}{^j_{ik}} X^k$). To raise indices, the inverse of the Rényi-metric \autoref{eq: renyi
metric} must be used, contributing an additional factor $\rho^{-1}$ as compared
to the Christoffel symbols of the first kind on the left. With \autoref{eq:
Amari Chentsov Tensor definition} and \autoref{eq: Levi-Civita connection
definition}, the following alternative formulations are possible,
\begin{align}
  \tensor[^{(\rho)}]{\nabla}{} 
  &= (1-\rho) \tensor[^{(e)}]{\nabla}{} + \rho \tensor[^{(m)}]{\nabla}{}
  =2 \rho \tensor[^{\text{LC}}]{\nabla}{} + (1-2\rho) \tensor[^{(e)}]{\nabla}{}= 2(1-\rho) \tensor[^{\text{LC}}]{\nabla}{} - (1-2\rho) \tensor[^{(m)}]{\nabla}{},
  \label{eq: rho connection in terms of e and C}\\  
  \tensor[^{(\rho^*)}]{\nabla}{}
  &= \rho \tensor[^{(e)}]{\nabla}{} + (1-\rho) \tensor[^{(m)}]{\nabla}{}
  = 2 \rho \tensor[^{\text{LC}}]{\nabla}{}+(1-2\rho)\tensor[^{(m)}]{\nabla}{} = 2 \left(1-\rho\right) \tensor[^{\text{LC}}]{\nabla}{} - (1-2\rho) \tensor[^{(e)}]{\nabla}{}\,.
  \label{eq: rho* connection in terms of e and C}
\end{align}

\subsection{Discussion}\label{sect: discussion renyi geometry}
This time, taking the limit of $\rho \to +1$ yields the Kullback-Leibler geometry,
\begin{equation}
\lim_{\rho \to 1} D_\rho
= D_{\text{KL}} \, ,
\quad
\lim_{\rho \to 1} \tensor[^{(\rho)}]{\nabla}{}
= \tensor[^{(e)}]{\nabla}{} + C
= \tensor[^{\text{KL}}]{\nabla}{}
\quad \text{and} \quad
\lim_{\rho \to 1} \tensor[^{(\rho^*)}]{\nabla}{}
= \tensor[^{(e)}]{\nabla}{}
= \tensor[^{\text{KL}^{*}}]{\nabla}{}. \label{eq: rho divergence connections for rho to 1}
\end{equation}
This asymptotic behavior differs from the Amari geometry ($\alpha \to -1$ but
$\rho \to +1$), due to the historic definitions of the divergences.

While the two geometries are
similar in the regard that they both admit a tensorial difference between the
connection $\nabla$ and its dual $\nabla^*$, which is given by a scaled version
of the Amari-Chentsov tensor $C$, they differ in two significant ways: On the
one hand, Amari's metric tensor is the Fisher information,
$\tensor[^{(\alpha)}]{g}{_{ij}} = \tensor[^{\text{F}}]{g}{_{ij}}$,
while the Rényi metric is additionally scaled by a conformal factor,
$\tensor[^{(\rho)}]{g}{_{ij}} = \rho \cdot \tensor[^{\text{F}}]{g}{_{ij}}$. Interestingly, the required positivity of the conformal factor is ensured by the definition of the Rényi divergence. With regards to the conformal invariance of Weyl-curvature, this conformal scaling is a promising object for future study.

On the other hand, while the 
$\alpha$-connections exhibit a symmetry in their parameter,
$\tensor[^{(\alpha^*)}]{\nabla}{} = \tensor[^{(-\alpha)}]{\nabla}{}$,
we lose this property for the case of
the Rényi connections $\tensor[^{(\rho^*)}]{\nabla}{} 
\neq \tensor[^{(-\rho)}]{\nabla}{}$.
(Still, 
$\tensor[^{(\rho^*)}]{\nabla}{}$ and $\tensor[^{(-\rho)}]{\nabla}{}$
fulfill dual metric compatibility with respect to ${}^{(\rho)}g$ since they are divergence-induced.)
Note
that the conformal factor does not affect the invariance under sufficient
statistics, since the Rényi metric is still proportional to the Fisher metric.

In fact, one may try to make the two geometries match by introducing a
coordinate transformation that alters the Rényi geometry,
\begin{equation}
\theta^{i'} = \sqrt{\rho}\, \theta^i
\rightarrow \tensor[^{(\rho)}]{g}{_{i'j'}} = \tensor[^{\text{F}}]{g}{_{i'j'}}\,.
\end{equation}
This undoes the conformal scaling and thus makes the two metrics match. However,
the connection coefficients
\begin{equation}
\tensor[^{(\rho)}]{\Gamma}{_{i'j'k'}} 
= \frac{1}{\rho^\frac{3}{2}}\tensor[^{(\rho)}]{\Gamma}{_{ijk}}
= \frac{1}{\sqrt{\rho}} \tensor[^{(e)}]{\Gamma}{_{ijk}} + \frac{1 - \rho}{\sqrt{\rho}} \tensor{C}{_{ijk}}
\end{equation}
do not take the same form as in the $\alpha$ geometry. Thus, while they share
some properties, the two geometries are genuinely different.

Lastly, consider the so-called Bhattacharyya distance~\cite{bhattacharyya_divergence_original_definition},
\begin{equation}
  D_{\text{B}}[\theta : \theta'] := D_{(\rho=\frac{1}{2})}[\theta : \theta'] = -2 \log \int \dd \mu(y) \,\sqrt{pp'}\,,
\end{equation}
a symmetric statistical divergence which emerges from Rényi (\autoref{eq: renyi divergence}) for the choice of $\rho = \frac{1}{2}$. Due to this divergence's symmetry, $D_{\text{B}}[\theta:\theta'] = D_{\text{B}}[\theta':\theta]$, its induced information geometry is \emph{non-dual}, i.e.
\begin{align}
  \tensor[^{\text{B}}]{g}{_{ij}} &= \frac{1}{2} \cdot \tensor[^{\text{F}}]{g}{_{ij}}, \\
  \tensor[^{\text{B}}]{\Gamma}{_{ijk}} 
  &= \frac{1}{2} \cdot \tensor[^{(e)}]{\Gamma}{_{ijk}} + \frac{1}{4} \cdot C_{ijk} = \tensor[^{\text{B}^*}]{\Gamma}{_{ijk}}  
  \; \leftrightarrow\; 
  \tensor[^{\text{B}}]{\nabla}{} = \tensor[^{\text{B}^*}]{\nabla}{} = \tensor[^{(e)}]{\nabla}{} + \frac{1}{2} \cdot C \,. 
\end{align}
and therefore we may identify
$\tensor[^{\text{B}}]{\nabla}{} 
= \tensor[^{\text{B}^*}]{\nabla}{}
= \tensor[^{\text{LC}}]{\nabla}{}$
by comparison with \autoref{eq: Levi-Civita connection definition}.
For an overview of the Bhattacharyya-geometry, please refer to \autoref{tab:
overview Bhattacharyya and KL geometry}. Clearly, the Bhattacharyya-geometry is
particular, in that its symmetry erases the usual duality of information geometry.
This is manifest in the self-dual connections $\tensor[^{\text{B}}]{\nabla}{} =
\tensor[^{\text{B}^*}]{\nabla}{} = \tensor[^{\text{LC}}]{\nabla}{}$ and so
forth.

\section{Dual Rényi-Laplace-Beltrami Operators}\label{sect: Dual Rényi-Laplace-Beltrami Operators} 
\subsection{Preliminaries: (Generalized) Vector Calculus in Information Geometry}
Generalized vector calculus on statistical manifolds is significant from both the geometric as well as the statistical angle. Consider again a manifold $\mathcal{M}$, a metric tensor $g$, as well as a connection $\nabla$. In differential geometry, the divergence of a vector field is defined as 
\begin{equation}
  \text{div}: \Gamma(T\mathcal{M}) \to \mathcal{F}(\mathcal{M}),\quad X \mapsto  \text{tr}\,\nabla X \;(=\nabla_{\partial_i} X^i) \,.\label{eq: diffgeo divergence definition}
\end{equation}
From a geometric standpoint, the divergence quantifies the rate at which the size of volume elements changes as they move along the flow of a vector field \cite{Udriste}. On a \emph{statistical} manifold, such changes in volume correspond to changes in the enclosed probability mass. Commonly, the LC-connection's divergence may be expressed as 
\begin{equation}
  \tensor[^{\text{LC}}]{\text{div}}{}X = \frac{1}{\sqrt{\det \tensor[]{g}{} }} \partial_j \left( \sqrt{\det \tensor[]{g}{}} X^j\right)\,,\label{eq: diffgeo LC divergence with metric tensor}
\end{equation}
in terms of a metric tensor $g$. Besides, the gradient of a smooth function,
\begin{equation}
  \text{grad}  : \mathcal{F}(\mathcal{M})\to \Gamma(T\mathcal{M}), \quad h \mapsto \,(\dd h)^{\sharp} \;\left(=g^{ij} \partial_j h\right)\,,
\end{equation} 
(with $\sharp$ the musical isomorphism) allows to define the generalized Laplacian or Laplace-Beltrami operator, 
\begin{equation}
  \Delta: \mathcal{F}(\mathcal{M}) \to \mathcal{F}(\mathcal{M}), \quad h \mapsto \text{div}\,\text{grad}\,h = \partial_i (\tensor[]{g}{^{ik}} \,\partial_k h) +  \tensor[]{\Gamma}{^j_{ij}}\;\tensor[]{g}{^{i\ell}} \,\partial_\ell h \,.\label{eq: diffgeo laplacian definition}
\end{equation}
Laplace-Beltrami operators are of great significance in statistics as they allow to evaluate the admissibility of statistical estimators, such as Bayes estimators, in the form of risk differences \cite{brown_admissible_1971,brown_heuristic_1979,hartiganmaximumlikelihoodprior,komaki_shrinkage_2001}. Consequently, it is our view that the parametrized Rényi-Laplacians derived in the following (in conjunction with the novel Rényi-priors, see Section \ref{sect: Rényi Volume Forms, Rényi-Priors}), may open alternative ways to judge the optimality of decision rules in statistics.

\subsection{The $\alpha$-geometry's Laplace-Beltrami Operator}\label{sect: The alpha-geometry's Laplace-Beltrami Operator}
Through \autoref{eq: diffgeo laplacian definition} with  $\tensor{g}{}=\tensor[^{\text{F}}]{g}{}$ and $\tensor[^{(\alpha)}]{\nabla}{}$, Calin and Udriste \cite{Udriste} establish the $\alpha$-divergence of a vector field, 
\begin{equation}
    \tensor[^{(\alpha)}]{\text{div}}{}X = \frac{1-\alpha}{2} \cdot \tensor[^{(m)}]{\text{div}}{}X + \frac{1+\alpha}{2} \tensor[^{(e)}]{\text{div}}{}X\,,\label{eq: alpha divergences of vector fields}
\end{equation}
as well as the the dual parametrized Laplace-Beltrami operator for Amari's $\alpha$-geometry ($\alpha$-Laplacian, $\alpha$-LB-operator) \cite{methodsamari,amari_information_2016},    
\begin{equation}
  \tensor[^{(\alpha)}]{\Delta}{} =  \tensor[^{\text{LC}}]{\Delta}{}  - \frac{\alpha}{2}\left(\tensor[^{(m)}]{\Delta}{} - \tensor[^{(e)}]{\Delta}{} \right) = \frac{1+\alpha}{2}\;\tensor[^{(e)}]{\Delta}{} + \frac{1-\alpha}{2}\; \tensor[^{(m)}]{\Delta}{}\,.  \label{eq: alpha Laplacian}
\end{equation}
Here, $\tensor[^{(e)}]{\text{div}}{}X,\tensor[^{(m)}]{\text{div}}{}X$ and $\tensor[^{(e)}]{\Delta}{}:=  \tensor[^{(e)}]{\text{div}}{}\,\tensor[^{\text{F}}]{\text{grad}}{}, \tensor[^{(m)}]{\Delta}{} :=  \tensor[^{(m)}]{\text{div}}{}\,\tensor[^{\text{F}}]{\text{grad}}{}$ are given by \autoref{eq: diffgeo divergence definition} and \autoref{eq: diffgeo laplacian definition} together with \autoref{eq: exponential connection coefficient defintion} and \autoref{eq: mixture connection coefficient and Amari Chentsov defintion}, respectively.

\subsection{The Rényi Dual Divergences and Rényi-Laplace-Beltrami-Operator}\label{sect: The Rényi Dual Divergences and Laplace-Beltrami Operator}
The following derivation follows Calin and Udriste's \cite{Udriste} derivation of \autoref{eq: alpha divergences of vector fields} and \autoref{eq: alpha Laplacian} (see also Amari and Nagaoka \cite{methodsamari}). In accordance with \autoref{eq: diffgeo divergence definition}, the dual Rényi divergences of a vector field $X \in \Gamma(T\mathcal{M})$ on a statistical manifold with $\tensor[^{(\rho)}]{\nabla}{},\tensor[^{(\rho^*)}]{\nabla}{}$ are denoted as
$\tensor[^{(\rho)}]{\text{div}}{}X = \tensor[^{(\rho)}]{\nabla}{}_{\partial_i} X^i$ and $\tensor[^{(\rho^*)}]{\text{div}}{}X = \tensor[^{(\rho^*)}]{\nabla}{}_{\partial_i} X^i $.

To begin with, we note that 
\begin{align}
   \tensor[^{\text{LC},(\rho)}]{\text{div}}{}X 
   \overset{\ref{eq: diffgeo LC divergence with metric tensor}}{=} \frac{1}{\sqrt{\det \tensor[^\rho]{g}{} }} \partial_j \left( \sqrt{\det  \tensor[^{\rho}]{g}{}}X^j\right) = \tensor[^{\text{LC},(\rho=1)}]{\text{div}}{}X \label{eq: div LC in terms of det g}\,, 
\end{align}
with $\tensor[^\rho]{g}{} = \rho\cdot \tensor[^{\text{F}}]{g}{}$ the metric tensor of the Rényi-geometry. Subsequently, \autoref{eq: rho connection in terms of e and C} and \autoref{eq: rho* connection in terms of e and C} lead to 
\begin{align}
   \tensor[^{(\rho)}]{\text{div}}{}X = \rho \cdot \tensor[^{(m)}]{\text{div}}{}X + (1-\rho) \tensor[^{(e)}]{\text{div}}{}X \label{eq: div rho in terms of div m and div e};\;
   \tensor[^{(\rho^*)}]{\text{div}}{}X = \rho \cdot \tensor[^{(e)}]{\text{div}}{}X + (1-\rho) \tensor[^{(m)}]{\text{div}}{}X \,.
\end{align}
Reparametrizing with $\rho(\alpha) = \frac{1-\alpha}{2}$ clarifies that the Rényi-divergences are equivalent to those of Amari's $\alpha$-geometry given in \autoref{eq: alpha divergences of vector fields}. The significance of this reparameterization will be the subject of a detailed discussion in Section \ref{sect: the Results volume forms}.

To find the dual RLB-operators, consider once again smooth function $h \in \mathcal{F}(\mathcal{M})$, the gradient of which can be expressed as  
\begin{align}
(\tensor[^{(\rho)}]{\text{grad}}{}\, h)^i &= \tensor[^{\rho}]{g}{^{ij}} \,\partial_j h =\rho^{-1}\cdot \tensor[^{\text{F}}]{g}{^{ij}} \,\partial_j h = \rho^{-1}\cdot (\tensor[^{\text{F}}]{\text{grad}}{}\, h)^i\,,\label{eq: rho grad definition}
\end{align}
for the Rényi-geometry. With $\tensor[^{(\rho)}]{\text{div}}{},\,\tensor[^{(\rho)}]{\text{grad}}{},\,$ in place, use of \autoref{eq: diffgeo laplacian definition} gives the RLB-operator 
\begin{align}
  \tensor[^{(\rho)}]{\Delta}{} h &:= \tensor[^{(\rho)}]{\text{div}}{}\,\tensor[^{(\rho)}]{\text{grad}}{}\, h \\
  &\overset{\eqref{eq: div rho in terms of div m and div e}}{=} \rho \cdot \tensor[^{(m)}]{\text{div}}{} (\tensor[^{(\rho)}]{\text{grad}}{} h) + (1-\rho) \tensor[^{(e)}]{\text{div}}{} (\tensor[^{(\rho)}]{\text{grad}}{} h) \\
  &\overset{\eqref{eq: rho grad definition}}{=}  \tensor[^{(m)}]{\Delta}{} h + \left(\rho^{-1}-1\right) \cdot \tensor[^{(e)}]{\Delta}{} h \,.\label{eq: rho Laplacian}
\end{align}
Note that in the limit of $\rho \to 1$, \autoref{eq: rho Laplacian} gives $\tensor[^{(\rho = 1)}]{\Delta}{}  =\tensor[^{(m)}]{\Delta}{}$. With this, one sees that $\tensor[^{(\rho)}]{\Delta}{} \neq \rho^{-1}\cdot \tensor[^{(\rho = 1)}]{\Delta}{} $. Conversely, the LC-connection's LB-operator scales with the inverse of the conformal factor, i.e. $\tensor[^{\text{LC},(\rho)}]{\Delta}{} := \tensor[^{\text{LC}}]{\text{div}}{} \tensor[^{(\rho)}]{\text{grad}}{} \overset{\eqref{eq: div LC in terms of det g}}{=} \rho^{-1} \cdot \tensor[^{\text{LC},(\rho = 1)}]{\Delta}{} $.

Direct comparison with \autoref{eq: alpha Laplacian} demonstrates that there is \emph{no way to reparameterize} the RLB-operator $\tensor[^{(\rho)}]{\Delta}{}$ to make it coincide with the $\alpha$-LB-operator $\tensor[^{(\alpha)}]{\Delta}{}$, i.e. 
\begin{equation}
  \tensor[^{(\rho)}]{\Delta}{}  =  \tensor[^{(m)}]{\Delta}{}  + \left(\rho^{-1} -1\right) \cdot \tensor[^{(e)}]{\Delta}{} \; \centernot\longleftrightarrow\; \tensor[^{(\alpha)}]{\Delta}{} = \frac{1-\alpha}{2}\; \tensor[^{(m)}]{\Delta}{} + \frac{1+\alpha}{2}\;\tensor[^{(e)}]{\Delta}{} \,,
\end{equation}
thus verifying Section \ref{sect: discussion renyi geometry}'s finding that the $\alpha$- and Rényi-geometry are genuinely different.

Looking forward, we wish to employ this novel parametrized Laplacian $\tensor[^{(\rho)}]{\Delta}{}$ to evaluate the admissibility of statistical estimators $\delta$ as well as predictive densities. As indicated by Hartigan \cite{hartiganmaximumlikelihoodprior}, Laplace-Beltrami-operators are commonly used to approximate \emph{risk differences} $R(\theta, \delta_\pi) - R(\theta, \delta_{\text{MLE}})$ of statistical estimators (with $\delta_\pi$ the Bayes estimator and $\delta_{\text{MLE}}$ the maximum likelihood estimator) as developed e.g. by Komaki \cite{komaki_shrinkage_2001}. Here, $\{\theta\}$ indicates a coordinate frame on the statistical manifold, while \mbox{$R(\theta, \delta) := \smallint \dd\mu(y)\, p_\theta(y) (\delta(y)-\theta)^2$} is the so-called \emph{risk} of the estimator in question (as defined in e.g. Brown \cite{brown_admissible_1971}, Komaki \cite{komaki_shrinkage_2001}). Moreover, in later works, Komaki \cite{Komaki_2006, komaki2015} uses the Laplacian $\Delta$ (for the LC-connection and the Fisher information metric) to quantify the risk difference between so-called predictive densities in statistical inference. 
In view of these promising applications, we plan to explore the effect that the $\rho$-scaling of the RLB-operators $\tensor[^{(\rho)}]{\Delta}{}$ (together with the Rényi-priors $\tensor[^{(\rho)}]{\text{cov}}{}, \tensor[^{(\rho^*)}]{\text{cov}}{}$ from Section \ref{sect: the Results volume forms}) has in deciding which estimator or predictive density is to be preferred over (i.e. has lower risk than) another.

\section{Rényi Volume Forms, Rényi-Priors}\label{sect: Rényi Volume Forms, Rényi-Priors}
\subsection{Preliminaries I: Covolumes as Priors}\label{sect: chapter_3_theory_priors}
Information geometry brought the geometric structure of statistical models to light, showing that parameter space is oftentimes non-Euclidean in nature. This important recognition comes with new challenges: $\dd^n \theta$ is no longer the correct integration measure on these curved manifolds \cite{amari_information_2016}.

To see why this is of key importance for statistics, consider e.g. the selection of a uniform prior in a Bayesian inference problem. In Bayesian inference, a prior belief (a probability distribution $\pi(\theta)$) shapes e.g. the evidence \mbox{$p(y) = \smallint \dd^n \theta\; \pi(\theta) \;p_\theta(y)$}, statistical estimators etc. \cite{amari_information_2016}. In cases where there is no a priori information about the parameters of interest, the most uninformative prior must be selected to ensure fairness. However, identifying the form of the uninformative (also called \enquote{uniform} in the literature) prior distribution over a non-Euclidean (possibly non-flat) statistical manifold is a highly non-trivial task. A prior which does \emph{not} correctly account for the statistical manifold's geometry contaminates the inference result, introducing unjustified bias in the form of e.g. artifacts of certain coordinate choices.

In Euclidean space, integration over an (improper) uniform prior is done by \mbox{$\smallint \dd^n \theta\; \pi(\theta) = \smallint \dd^n \theta \cdot 1$}. However, for a curved parameter space, the measure $\dd^n \theta$ does not assign equal probability mass to regions of equal volume, spoiling uniformity on the non-Euclidean statistical manifold \citep{Kass1996TheSO}. This issue must be fixed by a suitable prior choice, which acts as a weighting, a \emph{covolume}, that correctly restores the invariance of the integration measure, i.e. $\smallint \dd^n \theta \; \pi(\theta) := \smallint \omega$. Information geometry sometimes refers to such \enquote{covolume priors} as \emph{canonical}, since they incorporate the geometric structure of the statistical models at hand \citep{Matsuzoe2015InformationGO}. This Section derives the Rényi-geometry's priors, so as to correctly quantify the size of uncertainty regions for inference problems where the Rényi divergence is a favorable choice.

\subsection{Preliminaries II: Volume Forms in Differential Geometry}\label{sect: chapter_3_theory}
As before, consider an $n$-dimensional oriented Riemannian manifold $(\mathcal{M},g)$ \cite{Udriste}. A \emph{volume form} is a nowhere-vanishing $n$-form, which is given by  
\begin{equation}
  \omega = \text{cov}\; \dd x^1 \wedge \dots \wedge \dd x^n \;\in \Omega^n (\mathcal{M}) = \Gamma(\Lambda^n\,T^* \mathcal{M})\,,
\end{equation} 
with the local induced coframe $\{\dd x^i\}$ and $T^* \mathcal{M}$ the cotangent bundle \cite{nakahara2003geometry,lee2019introduction}. The coefficient function $\text{cov} \in C^\infty(\mathcal{M})$ is called the \emph{covolume}.

This work focuses on \emph{affine} volume forms, i.e. volume forms which are parallel with respect to a certain affine connection $\nabla$ \cite{alphaparallelpriors},   
  \begin{equation}
    \nabla \omega = 0\,,
\end{equation}
up to a constant factor, expressing the notion that volumes remain invariant under parallel transport. The canonical invariant volume form on a Riemannian manifold, the \emph{Riemannian volume form} is defined by 
    \begin{equation}
        \dd V := \sqrt{\det g}\; \dd x^1 \wedge \dots \wedge \dd x^n\,.  \label{eq: def Riemannian volume form in induced frame}
    \end{equation}
The Riemannian volume form's defining property is its parallelity with respect to the LC-connection, i.e. \mbox{$\tensor[^{\text{LC}}]{\nabla}{}\,\dd V = 0$}.

\subsection{Previous Studies on Geometric Priors}\label{sect: Existing Volume Forms (KL, Amari)}
Historically, it was Jeffreys who first introduced a canonical, \enquote{geometric} prior \mbox{$\pi_\text{J}(\theta) \sim \sqrt{\det  \tensor[^{\text{F}}]{g}{}(\theta)}$}, the so-called Jeffreys prior 
\cite{jeffreysprior101}. From a statistics standpoint, Jeffreys prior is the correct uniform prior distribution for an inference problem defined by the Fisher metric, due to its invariance under reparameterizations of the statistical model. From a geometric perspective, $\pi_J$ is the unique Riemannian covolume (\mbox{$\dd V_J = \pi_J \;\dd x^1\wedge \dots \wedge \dd x^n$}) on the potentially-curved statistical manifold. Thereby, Jeffreys volume form $\dd V_J $, which is parallel to the LC-connection, correctly defines integration over functions on the statistical manifold $\mathcal{M}$ \citep{jeffreysprior101, Nielsen_2020}.

Takeuchi and Amari's highly innovative work on geometric priors \cite{alphaparallelpriors} generalizes Jeffreys prior to a parametrized prior family induced by the dual, parametrized $\alpha$-connections $\tensor[^{(\alpha,\alpha^*)}]{\nabla}{}$. They define the $\alpha$-priors as the covolumes of the $\tensor[^{(\alpha)}]{\omega}{},\tensor[^{(\alpha^*)}]{\omega}{}$ volume forms defined by $\tensor[^{(\alpha)}]{\nabla}{}\tensor[^{(\alpha)}]{\omega}{} = 0$ $(\tensor[^{(\alpha^*)}]{\nabla}{} \tensor[^{(\alpha^*)}]{\omega}{} = 0)$. For the exponential family $\mathcal{P}_e$, they find \citep{alphaparallelpriors,Udriste} 
\begin{align}
    \pi_e^{(\alpha)} (\theta) := \tensor[^{(\alpha)}]{\text{cov}}{_e}\sim(\det \tensor[^{\text{F}}]{g}{})^{(1-\alpha)/2}\;\text{and}\;  \pi_e^{(\alpha^*)} (\theta) := \tensor[^{(\alpha^*)}]{\text{cov}}{_e}\sim(\det \tensor[^{\text{F}}]{g}{})^{(1+\alpha)/2}\,. \label{eq: alpha parallel covolumes}
\end{align}
Importantly, Takeuchi and Amari \cite{alphaparallelpriors} clarify that  $\tensor[^{(1)}]{\text{cov}}{_e}$ is the correct uniform prior distribution in the natural parameters $\theta$ (defined as the affine parameters with respect to $\tensor[^{(e)}]{\nabla}{}$, i.e. the local coordinate system in which the connection coefficients vanish, see \autoref{eq: exponential family}). Analogously, when one defines the expectation parameters $\{\eta^i\}$ as the affine coordinates of the $\tensor[^{(m)}]{\nabla}{}$ connection on $\mathcal{P}_e$, then $\tensor[^{(1^*)}]{\text{cov}}{_e}$ is the uniform prior distribution in $\{\eta^i\}$ \cite{alphaparallelpriors} (see \autoref{eq: mixture family}). Accordingly, the covolume of the volume element parallel to $\tensor[^{(\alpha)}]{\nabla}{}$ is distributed uniformly in the $\alpha$-affine coordinates\footnote{\emph{\enquote{If there exist affine coordinates with respect to $\nabla$, a density on $\mathcal{M}$ induced by $\kappa$ is uniform with respect to those affine coordinates.}} (Quote from Amari and Takeuchi \cite{alphaparallelpriors}, p. 1015)}. Besides, the $\alpha$-priors recover Jeffreys prior for the choice of $\alpha = 0$ \cite{alphaparallelpriors}.

In the limit of $\alpha \to 1$, the volume forms parallel to $\tensor[^{(\alpha=1), (\alpha^*=1^*)}]{\nabla}{} = \tensor[^{(e),(m)}]{\nabla}{}$ are
\begin{align}
  \tensor[^{(1)}]{\text{cov}}{_e} \;\dd x^1\wedge \dots \wedge \dd x^n \propto \dd x^1\wedge \dots \wedge \dd x^n  \;&\text{parallel to}\; \tensor[^{(e)}]{\nabla}{}\,,\\
  \tensor[^{(1^*)}]{\text{cov}}{_e} \;\dd x^1\wedge \dots \wedge \dd x^n \propto (\det \tensor[^{\text{F}}]{g}{}) \;\dd x^1\wedge \dots \wedge \dd x^n\; &\text{parallel to}\; \tensor[^{(m)}]{\nabla}{}\,.\label{eq: alpha covolumes for alpha to 1}
\end{align} 
In the limit of $\alpha \to 1$, the $\alpha$-divergence recovers the Kullback-Leibler divergence. Accordingly, the limit of $\alpha \to 1$ gives the covolumes of the KL-divergence. This result is quickly confirmed by directly deriving the KL-covolumes for its induced geometry. As expected, the covolumes for the KL-divergence-induced connections $\tensor[^{\text{KL}, {\text{KL}}^*}]{\nabla}{} = \tensor[^{(m),(e)}]{\nabla}{}$ are the uniform prior densities in the natural (expectation) parameters, respectively \citep{hartiganmaximumlikelihoodprior,alphaparallelpriors}, i.e. up to a constant 
\begin{align}
  \pi_e^{\text{KL}} (\theta) := \tensor[^{\text{KL}}]{\text{cov}}{_e} \sim (\det \tensor[^{\text{F}}]{g}{}) \quad \text{and}\quad 
  \pi_e^{\text{KL}^*} (\theta):= \tensor[^{\text{KL}^*}]{\text{cov}}{_e} \sim 1\,,
\end{align}
with 
\begin{align}
  \tensor[^{\text{KL}}]{\text{cov}}{_e} \; \dd x^1\wedge \dots \wedge \dd x^n \;\propto \;(\det \tensor[^{\text{F}}]{g}{})  \;\dd x^1\wedge \dots \wedge \dd x^n  &\;\text{parallel to}\;\tensor[^{(m)}]{\nabla}{}\,,\\
  \tensor[^{{\text{KL}}^*}]{\text{cov}}{_e} \;\dd x^1\wedge \dots \wedge \dd x^n\; \propto\;\dd x^1\wedge \dots \wedge \dd x^n &\;\text{parallel to}\;  \tensor[^{(e)}]{\nabla}{}\,, \label{eq: parallelity of KL volume forms}
\end{align} 
which is in agreement with \autoref{eq: alpha covolumes for alpha to 1}, as required. Various studies have reviewed canonical geometric priors and their applications, see e.g. \citep{jeffreysprior101, alphaparallelpriors, Matsuzoe2015InformationGO, weylprior}.

\subsection{The Rényi-Priors and Hartigan's Prior Family}\label{sect: the Results volume forms}
This Section derives the volume forms that are parallel to the dual Rényi-induced connections given in Section \ref{sect: Rényi's rho geometry}, as well as the corresponding covolumes which we call the \emph{Rényi-priors}. The technical computation of these volume forms is based on the discussion of Takeuchi and Amari's $\alpha$-priors (see \autoref{eq: alpha parallel covolumes}) presented in Calin and Udriste \citep{Udriste}. For a detailed derivation, the reader is referred to \ref{sect: Derivation of the Rényi Volume Forms}. The link of these Rényi-priors to a prior family proposed by Hartigan \cite{hartiganmaximumlikelihoodprior} is elucidated.

To begin with, consider the Rényi-volume-forms,  
\begin{align}
  \tensor[^{(\rho)}]{\omega}{} =  \tensor[^{(\rho)}]{\text{cov}}{}\; \dd x^1 \wedge \dots \wedge \dd x^n\,,\quad  
  \tensor[^{(\rho^*)}]{\omega}{} = \tensor[^{(\rho^*)}]{\text{cov}}{}\; \dd x^1 \wedge \dots \wedge \dd x^n\,, \label{eq: def Rényi volume forms}
\end{align}
defined by the partial differential equations 
  \begin{equation}
    \tensor[^{(\rho)}]{\nabla}{} \tensor[^{(\rho)}]{\omega}{}  = 0\,, \quad \tensor[^{(\rho^*)}]{\nabla}{} \tensor[^{(\rho^*)}]{\omega}{}  = 0\,. \label{eq: general definitions of Rényi volume forms}
  \end{equation}
On the oriented statistical manifold $(\mathcal{M},\tensor[^\rho]{g}{})$ with $\tensor[^\rho]{g}{} = \rho \cdot \tensor[^{\text{F}}]{g}{}$, an $n$-form $\tensor[^{(\rho)}]{\omega}{} \in \Omega^n(\mathcal{M})$ can be written as \cite{Udriste}
\begin{equation}
    \tensor[^{(\rho)}]{\omega}{} := f^\rho \cdot \dd V_\rho,\label{eq: Rényi volume form in terms of f_rho}
\end{equation} 
with $f_\rho \in \mathcal{F}(\mathcal{M})$ a smooth nowhere-vanishing function and $\dd V_\rho$ 
the Riemannian volume form,
\begin{equation}
    \dd V_\rho := \rho^{n/2} \cdot \sqrt{\det \tensor[^{\text{F}}]{g}{}} \;\dd x^1 \wedge \dots \wedge \dd x^n\,, 
\end{equation} 
which is parallel to $\tensor[^{(\rho), \text{LC}}]{\nabla}{} $.

We first consider the exponential family $\mathcal{P}_e$'s statistical manifold. The dual Rényi-volume forms for $\mathcal{P}_e$ are referred to as $\tensor[^{(\rho)}]{\omega}{}_e := f_e^\rho \, \dd V_\rho\;\text{and}\;  \tensor[^{(\rho^*)}]{\omega}{}_e := f_e^{\rho^*} \, \dd V_\rho$, respectively. The investigation begins with the $\tensor[^{(\rho^*)}]{\nabla}{}$-parallel volume form. Combining \autoref{eq: general definitions of Rényi volume forms} and \autoref{eq: Rényi volume form in terms of f_rho} leads to 
\begin{equation}
  -\partial_i (\log f^{\rho^*}_e)\,\dd V_\rho = \tensor[^{(\rho^*)}]{\nabla}{_{\partial_i}} \dd V_{\rho} \,.
\end{equation}
\ref{sect: Derivation of the Rényi Volume Forms} details how the computation of the covariant $\rho^*$-derivative of $\dd V_\rho$ is combined with the exponential family's flatness in $\tensor[^{(e)}]{\nabla}{}$ to reach 
\begin{equation}
  0 = \partial_i \log \left(f_e^{\rho^*} (\det \tensor[^{\text{F}}]{g}{})^{+ \left(\rho-\frac{1}{2}\right)}\right) \quad 
  \leftrightarrow \quad f_e^{\rho^*} = (\det \tensor[^{\text{F}}]{g}{})^{- \left(\rho-\frac{1}{2}\right)}\,,
\end{equation}
up to a constant factor. Finally, one arrives at the $\rho^*$-parallel volume form, 
\begin{equation}
    \tensor[^{(\rho^*)}]{\omega}{}_e = \rho^{n/2}\cdot (\det \tensor[^{\text{F}}]{g}{})^{1-\rho} \;\dd x^1\wedge \dots \wedge \dd x^n\,. \label{eq: dual volume form exponential family}
\end{equation}  
In analogy, its dual, the $\rho$-parallel volume form can derived 
\begin{equation}
 \tensor[^{(\rho)}]{\omega}{}_e = \rho^{n/2}\cdot (\det \tensor[^{\text{F}}]{g}{})^{+\rho} \;\dd x^1\wedge \dots \wedge \dd x^n\,.\label{eq: volume form exponential family}
\end{equation}  
It follows that the $\rho^*,\rho$-parallel covolumes, the Rényi-priors, on $\mathcal{P}_e$ turn out to be  
\begin{equation}
  \tensor[^{(\rho^*)}]{\text{cov}}{_e}=\rho^{n/2}\cdot(\det \tensor[^{\text{F}}]{g}{})^{1-\rho} \quad \text{and} \quad \tensor[^{(\rho)}]{\text{cov}}{_e}=\rho^{n/2}\cdot(\det \tensor[^{\text{F}}]{g}{})^{+\rho} \,.\label{eq: Renyi priors exponential family}
\end{equation}
Analogously, one can examine dual volume forms for the mixture family $\mathcal{P}_m$'s statistical manifold ($\tensor[^{(\rho)}]{\omega}{}_m, \tensor[^{(\rho^*)}]{\omega}{}_m$) to find $\tensor[^{(\rho^*)}]{\omega}{}_m, \tensor[^{(\rho)}]{\omega}{}_m$ with their covolumes
\begin{equation}
  \tensor[^{(\rho^*)}]{\text{cov}}{_m}=\rho^{n/2}\cdot(\det \tensor[^{\text{F}}]{g}{})^{+\rho} \quad \text{and} \quad \tensor[^{(\rho)}]{\text{cov}}{_m}=\rho^{n/2}\cdot(\det \tensor[^{\text{F}}]{g}{})^{1-\rho} \,. \label{eq: Renyi priors mixture family}
\end{equation}
From these results, several interesting findings emerge:
\begin{enumerate}
  \item It is evident that the Rényi-priors coincide with Jeffreys prior for $\rho = \frac{1}{2}$. Also, the limit of $\rho\to\frac{1}{2}$ in \autoref{eq: Renyi priors exponential family} (and \autoref{eq: Renyi priors mixture family}) is distinct in that this is the only choice for which the volume forms are symmetric $\tensor[^{(\rho=(1/2))}]{\omega}{}_{e,m} = \tensor[^{(\rho^*=(1/2)^*)}]{\omega}{}_{e,m}$, 
\begin{align}
  \tensor[^{(\rho=(1/2))}]{\omega}{}_{e,m} = \tensor[^{(\rho^*=(1/2)^*)}]{\omega}{}_{e,m} \sim \sqrt{\det \tensor[^{\text{F}}]{g}{}} \;\dd x^1 \wedge \dots \wedge \dd x^n \,.
\end{align}
Of course, this is in accordance with expectations, as $\rho = \frac{1}{2}$ makes the Rényi divergence coincide with the \emph{symmetric} Bhattacharyya divergence $D_{\text{B}}$, see Section \ref{sect: Rényi's rho geometry}. It is straightforward to show that the Jeffreys prior $\pi_J$ is the canonical uniform prior in the (self-dual) Bhattacharyya-connection's affine coordinates. 
\item The Rényi-priors recover the KL-divergence-induced volume forms for $\rho = 1$ (which are also reached for $\alpha = -1$), reproducing the $\rho\to 1$ ($\alpha\to -1$) limit at the level of divergences \citep{van_erven_renyi_2014}. In fact, comparison of \autoref{eq: Renyi priors exponential family} in the limit of $\rho \to 1$ with \autoref{eq: alpha covolumes for alpha to 1} allows for a consistency check \cite{alphaparallelpriors, Udriste}: as required, 
  $\tensor[^{(\rho^* = 1^*)}]{\text{cov}}{_e} = 1$ is uniform in the affine coordinates of $\tensor[^{(\rho^* = 1^*)}]{\nabla}{}  = \tensor[^{(e)}]{\nabla}{}$ and $\tensor[^{(\rho = 1)}]{\text{cov}}{_e} = (\det \tensor[^{\text{F}}]{g}{})$ is uniform in the affine coordinates of $\tensor[^{(\rho = 1)}]{\nabla}{}  = \tensor[^{(m)}]{\nabla}{}$, supporting the validity of the Rényi-priors. 
\item The Rényi-priors demonstrate the duality of the exponential and mixture family at the level of volume forms, as $ \tensor[^{(\rho^*)}]{\omega}{}_e = \tensor[^{(\rho)}]{\omega}{}_m$ and $\tensor[^{(\rho)}]{\omega}{}_e = \tensor[^{(\rho^*)}]{\omega}{}_m$. This is to be expected, as it replicates the duality of the $\tensor[^{(e)}]{\nabla}{}$- and $\tensor[^{(m)}]{\nabla}{}$-connections \cite{manyfaces}.
\item The fact that the $\rho$-covolumes in \autoref{eq: Renyi priors exponential family} and \autoref{eq: Renyi priors mixture family} differ from the $\alpha$-covolumes (\autoref{eq: alpha parallel covolumes}) by a factor of $\rho^{n/2}$ is expected, since conformal changes of the metric cause a scaling with the $(n/ 2)$-th power of the conformal factor in the volume form. This is however insignificant as volume forms are only defined up to a constant factor. 
\item Comparison of the Rényi-priors with \autoref{eq: alpha parallel covolumes} shows that, up to a power of the conformal factor, the dual volume forms induced by the Rényi divergence can be translated to those derived from $\alpha$-divergences (\autoref{eq: alpha parallel covolumes}) by the parameter choice $\rho = \frac{(1-\alpha)}{2}$. For the covolumes in \autoref{eq: Renyi priors exponential family}, one finds for instance 
\begin{align}
\tensor[^{(\rho)}]{\text{cov}}{_e} = \tensor[^{(\rho^*)}]{\text{cov}}{_m} &\sim (\det \tensor[^{\text{F}}]{g}{})^{+\rho} \; \overset{\rho = (1-\alpha)/2}{\longrightarrow} \;   \tensor[^{(\alpha)}]{\text{cov}}{_e} = (\det \tensor[^{\text{F}}]{g}{})^{(1-\alpha)/2} \,,  \label{eq: cov rho to cov alpha} \\
  \tensor[^{(\rho^*)}]{\text{cov}}{_e} =  \tensor[^{(\rho)}]{\text{cov}}{_m} &\sim  (\det \tensor[^{\text{F}}]{g}{})^{1-\rho} \; \overset{\rho = (1-\alpha)/2}{\longrightarrow} \; \tensor[^{(\alpha^*)}]{\text{cov}}{_e} \; = (\det \tensor[^{\text{F}}]{g}{})^{(1+\alpha)/2}\,. \label{eq: dual cov rho to dual cov alpha}
\end{align}
This result is non-trivial in light of the fact that the Rényi divergence's logarithmic structure does not emerge from the $\alpha$-divergence by merely choosing this alternative parametrization. Furthermore, Sections \ref{sect: Rényi Metric and Connections} and \ref{sect: The Rényi Dual Divergences and Laplace-Beltrami Operator} clarifies that the geometries are not the same. 
\item Crucially, Amari and Takeuchi \cite{alphaparallelpriors} recognize the reparameterization the $\alpha$-priors, $\alpha' = \frac{(1-\alpha)}{2}$, as a special choice for which certain asymptotic equivalences of statistical estimators are reached. The projected Bayes estimator asymptotically coincides with the bias-corrected maximum likelihood estimator in this reparameterization of the $\alpha$-priors \cite{alphaparallelpriors}. Technically, this reparameterization is an artifact of $\tensor[^{(\alpha)}]{\Gamma}{_{ijk}}$'s scaling with the Amari-Chentsov tensor $C_{ijk}$ in \autoref{eq: scaling of alpha connection with C}. Moreover, Takeuchi and Amari \cite{alphaparallelpriors} state that the $\alpha$-priors coincide with a parametrized prior family introduced by Hartigan for this choice \cite{hartiganmaximumlikelihoodprior, alphaparallelpriors}, thus giving geometric validation to Hartigan's statistics result. However, this is presented as a phenomenological result, without an a priori geometric reason as to why $\rho = \frac{(1-\alpha)}{2}$ \emph{should} be of particular significance \cite{alphaparallelpriors}. 

We thus come to our main result. \autoref{eq: cov rho to cov alpha} and \autoref{eq: dual cov rho to dual cov alpha} show that the Rényi-priors offer a geometric rationale: the choice of $\rho = \frac{(1-\alpha)}{2}$ gives the covolumes of the Rényi-geometry for $\mathcal{P}_e$ ($\mathcal{P}_m$). To put it differently, Hartigan's prior family coincides with the family of priors induced by the Rényi-geometry. The remainder of this Section is reserved for a derivation of this result. 
\end{enumerate}

In his 1964 and 1965 \cite{hartiganinvariantpriordistributions, hartigan1965} works on asymptotically unbiased prior distributions in statistics, Hartigan introduced a parametrized family of priors $\{\pi_{\alpha_H}\}$. Notably, Hartigan \cite{hartiganinvariantpriordistributions, hartigan1965} motivates this prior family solely from a statistics point of view, not a geometric one. In his later work \cite{hartiganmaximumlikelihoodprior}, Hartigan gives the following definition for this parametrized family of invariant priors,
\begin{equation}
    \partial_i \log \pi_{\alpha_H}:= (\tensor[^{\text{F}}]{g}{^{-1}})^{jk} \;\mathbb{P}_\theta \left[\alpha_H \cdot \partial_i \log  p_\theta \;\partial_j \log p_\theta\; \partial_k \log p_\theta  + \partial_i \partial_j \log  p_\theta \;\partial_k \log  p_\theta \right]\,.\label{eq: hartigans prior family definition}
\end{equation} 
We denote Hartigan's parameter by $\alpha_H$ to distinguish the $\alpha_H$-priors from Takeuchi and Amari's $\alpha$-priors. Moreover, we assume Hartigan's notation $\mathbb{P}_\theta[f]$ for the conditional expected value of a function over data space, i.e. $\mathbb{P}_\theta[f] := \smallint \dd \mu(y) \;\mathcal{L}(y|\theta) \, f(y)$. The following derivation of this result is not limited to $\mathcal{P}_e$ or $\mathcal{P}_m$, but holds true for statistical manifolds constructed from any suitable family of probability distributions.

Evaluation the parallelity condition, $0 = (\tensor[^{(\rho)}]{\nabla}{}\tensor[^{(\rho)}]{\omega}{})$, for an affine Rényi volume form $\tensor[^{(\rho)}]{\omega}{}$ in a local induced frame $\{\partial_i\}$ leads to 
\begin{align}
    \partial_i(\log \tensor[^{(\rho)}]{\text{cov}}{}) &= \sum_k \dd x^1 \wedge\dots\wedge \dd x^n (\partial_1,\dots,\tensor[^{(\rho)}]{\Gamma}{^j_{ki}} \partial_j, \dots,\partial_n) =  \tensor[^{(\rho)}]{\Gamma}{^j_{ji}} \,.\label{eq: hartigan and rényi}
\end{align}
Use of \autoref{eq: rho connection in terms of e and C} gives $\partial_i(\log \tensor[^{(\rho)}]{\text{cov}}{}) = \rho \cdot C_{ijk} + \tensor[^{(e)}]{\Gamma}{^j_{ij}}$.
Finally, we insert \autoref{eq: mixture connection coefficient and Amari Chentsov defintion} and \autoref{eq: exponential connection coefficient defintion} to find 
\begin{align}
  \partial_i(\log \tensor[^{(\rho)}]{\text{cov}}{}) = (\tensor[^{\text{F}}]{g}{^{-1}})^{jk} \;\mathbb{P}_\theta \left[\rho \cdot \partial_i \log  p_\theta \;\partial_j \log p_\theta\; \partial_k \log p_\theta  + \partial_i \partial_j \log  p_\theta \;\partial_k \log  p_\theta \right] \,, \label{eq: renyi covolume in hartigan notation}
\end{align}
Comparison of \autoref{eq: renyi covolume in hartigan notation} and \autoref{eq: hartigans prior family definition} immediately shows that $\alpha_H = \rho$, i.e. Hartigan's prior family with its original parametrization, is precisely the prior family of the Rényi-geometry. This result sustains the findings of Takeuchi and Amari \cite{alphaparallelpriors} whilst offering a new explanation as to why the reparameterization $\alpha' = \frac{(1-\alpha)}{2}$ holds particular geometric and statistical meaning: this choice recovers Rényi's uniform priors.

\section{Conclusions and Outlook} \label{sect: Conclusions and Outlook}
Building on the results by Eguchi \cite{eguchi1983divergenceinduced}, Jeffreys \cite{jeffreysprior101}, as well as Takeuchi and Amari \cite{alphaparallelpriors,amari_information_2016}, this work explored the information geometry induced by the statistical Rényi divergence. Particular emphasis was placed on the role of certain geometric quantities, such as dual covolumes, in the context of statistical inference. Our main conclusions and findings can be stated as follows: 
\begin{enumerate}
  \item Both, the metric tensor as well as the dual connections were computed from the Rényi-divergence. The Rényi-geometry's metric tensor was found to differ from the Fisher information metric by a conformal factor $\rho$. In spite of the fact that the dual Rényi-connections are not symmetric in their parameters in the sense of $\tensor[^{(\rho^*)}]{\nabla}{} 
  \neq \tensor[^{(-\rho)}]{\nabla}{}$, dual metric compatibility holds, as is requisite for an operable information geometry. In general, the dual Rényi-geometry is genuinely different from Amari's $\alpha$-geometry, as indicated by the failure of any general coordinate transformation to make the parametrized $\rho$-scaling and $\alpha$-scaling, respectively, coincide. As required, these geometries coincide in the limit of $(\rho \to 1), \,(\alpha \to -1)$, respectively. 
  \item Dual parametrized Laplace-Beltrami-operators were calculated from the Rényi-geometry's dual connections. The structure of these operators exhibits a $\rho$-scaling which \emph{cannot} be made to coincide with that of Amari's $\alpha$-LB-operator by way of a reparameterization $\rho(\alpha)$, supporting the recognition that these information geometries are inherently different in structure. 
  \item The dual Rényi-priors were derived as covolumes of the dual volume forms, which are parallel to the individual dual Rényi-connections. The results confirm the duality of the exponential and mixture family \cite{manyfaces} at the level of volume forms, since $ \tensor[^{(\rho^*)}]{\omega}{}_e = \tensor[^{(\rho)}]{\omega}{}_m$ and $\tensor[^{(\rho)}]{\omega}{}_e = \tensor[^{(\rho^*)}]{\omega}{}_m$. As required, the Rényi-priors recover the self-dual Jeffreys prior for $\rho \to (1/2)$ and coincide with Amari and Takeuchi's $(\alpha = 1)$-priors for $(\rho = 1)$. Up to a power of the conformal factor, the reparameterization $\rho(\alpha) = (1+\alpha)/2$ translates the canonical Rényi-priors to the $\alpha$-priors. Importantly, Amari and Takeuchi note that the choice $\alpha \to (1+\alpha)/2$ is distinct, as it leads to certain unique asymptotic behaviours of statistical estimators, all while also recovering Hartigan's prior family \cite{alphaparallelpriors,hartiganmaximumlikelihoodprior}. We find that for this choice, the $\alpha$-priors become the Rényi-priors. To put it differently, Hartigan's prior family with its original $\alpha_H$-parametrization precisely coincides with the Rényi-priors (i.e. $\rho = \alpha_H$). 
\end{enumerate}
The following table presents the essential geometric quantities derived for the Rényi-geometry in this work, in direct comparison with those of some well-established information geometries, namely the Bhattacharyya-geometry, the KL-geometry, as well as Amari's $\alpha$-geometry.  
\begin{table}[H]
  \centering
    \small 
    \setlength{\tabcolsep}{9pt} 
    \renewcommand{\arraystretch}{1.5} 
    \begin{tabular}{c c c}
    \hline
    & Bhattacharyya Divergence & Kullback-Leibler Divergence \\
    \hline
    Divergence & $D_{\text{B}}[\theta : \theta'] = - 2\log \int \dd \mu(y)\,\sqrt{p \,p'}$ & $D_{\text{KL}}[\theta : \theta'] = \int \dd \mu(y) \, p \log \frac{p}{p'} $ \\
    Metric & $\tensor[^{\text{B}}]{g}{_{ij}} =  \frac{1}{2}\cdot\tensor[^{\text{F}}]{g}{_{ij}}$ & $\tensor[^{\text{KL}}]{g}{_{ij}} = \tensor[^{\text{F}}]{g}{_{ij}}$\\
    Conn.& $\tensor[^{\text{B}}]{\nabla}{}= \tensor[^{\text{B}^*}]{\nabla}{}=  \tensor[^{\text{LC}}]{\nabla}{} $ & $\tensor[^{\text{KL}}]{\nabla}{}= \tensor[^{(m)}]{\nabla}{}$ \\
    Dual Conn. & $\tensor[^{\text{B}^*}]{\nabla}{}= \tensor[^{\text{B}}]{\nabla}{}= \tensor[^{\text{LC}}]{\nabla}{}$ & $\tensor[^{\text{KL}^*}]{\nabla}{} = \tensor[^{(e)}]{\nabla}{}$ \\
    $\mathcal{P}_e$ Cov. & $\tensor[^{\text{B}}]{\text{cov}}{_e} = \tensor[^{\text{LC}}]{\text{cov}}{_e} \propto \sqrt{\det \tensor[^{\text{F}}]{g}{}}$ & $\tensor[^{\text{KL}}]{\text{cov}}{_e}\propto(\det \tensor[^{\text{F}}]{g}{})$ \\
    $\mathcal{P}_e$ Dual Cov.& $\tensor[^{\text{B}^*}]{\text{cov}}{_e} = \tensor[^{\text{LC}}]{\text{cov}}{_e} \propto \sqrt{\det \tensor[^{\text{F}}]{g}{}}$ & $\tensor[^{\text{KL}^*}]{\text{cov}}{_e} \propto1$ \\
    $\mathcal{P}_m$ Cov. & $\tensor[^{\text{B}}]{\text{cov}}{_m} = \tensor[^{\text{LC}}]{\text{cov}}{_m} \propto\sqrt{\det \tensor[^{\text{F}}]{g}{}}$ & $\tensor[^{\text{KL}}]{\text{cov}}{_m} \propto 1$ \\
    $\mathcal{P}_m$ Dual Cov.& $\tensor[^{\text{B}^*}]{\text{cov}}{_m} = \tensor[^{\text{LC}}]{\text{cov}}{_e} \propto\sqrt{\det \tensor[^{\text{F}}]{g}{}}$ & $\tensor[^{\text{KL}^*}]{\text{cov}}{_m} \propto(\det \tensor[^{\text{F}}]{g}{})$ \\
    \hline
    \end{tabular}
    \caption{Overview of the Bhattacharyya-geometry alongside the KL-geometry.}
    \label{tab: overview Bhattacharyya and KL geometry}
  \end{table}

    \begin{table}[H]
      \centering
        \small 
        \setlength{\tabcolsep}{9pt} 
        \renewcommand{\arraystretch}{1.5}
        \begin{tabular}{c c c}
        \hline
        & Amari $\alpha$-Divergence & Rényi $\rho$-Divergence \\
        \hline
        Divergence &  $D_\alpha[\theta : \theta'] = \frac{4}{1-\alpha^2} \left(
1 - \int \dd \mu(y) p^\frac{1-\alpha}{2} (p')^\frac{1 + \alpha}{2}
\right) $ &  $D_{\rho}[\theta : \theta'] = \frac{1}{\rho - 1} \log \int \dd\mu(y) p^\rho (p')^{1-\rho} $ \\
        Metric & $\tensor[^{(\alpha)}]{g}{_{ij}} = \tensor[^{\text{F}}]{g}{_{ij}}$ & $\tensor[^{(\rho)}]{g}{_{ij}} = \rho\cdot\tensor[^{\text{F}}]{g}{_{ij}}$ \\
        Conn.& $\tensor[^{(\alpha)}]{\nabla}{}  = \tensor[^{(e)}]{\nabla}{} + \frac{(1-\alpha)}{2} C$ & $\tensor[^{(\rho)}]{\nabla}{}  = \tensor[^{(e)}]{\nabla}{} + \rho\, C$\\
        Dual Conn. & $\tensor[^{(\alpha^*)}]{\nabla}{}  = \tensor[^{(e)}]{\nabla}{} + \frac{1+\alpha}{2} C$ & $\tensor[^{(\rho^*)}]{\nabla}{}   =  \tensor[^{(e)}]{\nabla}{} + (1-\rho)\, C$\\
        $\mathcal{P}_e$ Cov. & $\tensor[^{(\alpha)}]{\text{cov}}{_e}=(\det \tensor[^{\text{F}}]{g}{})^{(1-\alpha)/2}$ & $\tensor[^{(\rho)}]{\text{cov}}{_e}=\rho^{n/2}\cdot(\det \tensor[^{\text{F}}]{g}{})^\rho$ \\
        $\mathcal{P}_e$ Dual Cov.& $\tensor[^{(\alpha^*)}]{\text{cov}}{_e}=(\det \tensor[^{\text{F}}]{g}{})^{(1+\alpha)/2}$ & $\tensor[^{(\rho^*)}]{\text{cov}}{_e}=\rho^{n/2}\cdot(\det \tensor[^{\text{F}}]{g}{})^{1-\rho}$  \\
        $\mathcal{P}_m$ Cov. & $\tensor[^{(\alpha)}]{\text{cov}}{_m}=(\det \tensor[^{\text{F}}]{g}{})^{(1+\alpha)/2}$ & $\tensor[^{(\rho)}]{\text{cov}}{_m}=\rho^{n/2}\cdot(\det \tensor[^{\text{F}}]{g}{})^{1-\rho}$ \\
        $\mathcal{P}_m$ Dual Cov.& $\tensor[^{(\alpha^*)}]{\text{cov}}{_m}=(\det \tensor[^{\text{F}}]{g}{})^{(1-\alpha)/2}$ & $\tensor[^{(\rho^*)}]{\text{cov}}{_m}=\rho^{n/2}\cdot(\det \tensor[^{\text{F}}]{g}{})^{\rho}$ \\
        \hline
        \end{tabular}
        \caption{Overview of Amari's $\alpha$-geometry \cite{amari_information_2016,alphaparallelpriors} alongside the Rényi-geometry.}
        \label{tab: geometric quantities alpha renyi}
        \end{table}
In the future, we plan to explore the integration measure defined by Rényi's volume forms, i.e. for an open set $U\subseteq \mathcal{M}$ on the statistical manifold, define integration as  
\begin{equation}
  \mu_{\tensor[^{(\rho),(\rho^*)}]{\omega}{}}(U) = \int_U \tensor[^{(\rho),(\rho^*)}]{\omega}{}\,.
\end{equation}
Furthermore, future research is needed to determine how the extremality of the family of Gaussian probability distributions arises in the context of volume integration on the statistical manifold. In this effort, we plan to examine the link between isoprobability surfaces and enclosed probability mass for the geometry of different statistical families. As hinted at in Section \ref{sect: Dual Rényi-Laplace-Beltrami Operators}, more work is needed to examine the dual RLB-operator's effects in judging the optimality of statistical estimators and predictive densities. In addition, we plan to explicate the link between RLB-operators and Rényi-covolumes. Further topics of interest for the Rényi-geometry include Lie derivatives and Killing fields for quantifying statistical degeneracies as isometries, particularly since the conventional definition of the Lie derivative relies on metric compatibility of the covariant derivative.

\vspace{1cm} 

\section*{Author Contributions} Foundations of Rényi-Geometry (Metric Tensor, Connections): HvC, Advanced Geometrical Aspects  (Rényi-Laplacian and Rényi-priors): RMK, Overall Conceptualization and Supervision: BMS.

\section*{Funding} This work was supported by the Deutsche Forschungsgemeinschaft (DFG, German
Research Foundation) under Germany's Excellence Strategy EXC 2181/1 - 390900948
(the Heidelberg STRUCTURES Excellence Cluster). RMK acknowledges funding of the
Stiftung der Deutschen Wirtschaft (Foundation of German Business) with funds
from the Begabtenförderung of the BMBF (Federal Ministry of Education and
Research's scholarship programme for gifted students). HvC
 is supported by the Konrad Zuse School of Excellence in Learning
and Intelligent Systems (ELIZA) through the DAAD programme Konrad Zuse Schools
of Excellence in Artificial Intelligence, sponsored by the Federal Ministry of
Education and Research.

\section*{Acknowledgments}
The authors wish to thank Benedikt Schosser for helpful discussions.

\section*{Conflicts of Interest}
The authors declare no conflicts of interest.

\appendix
\section{Detailed Derivation of the Rényi Volume Forms}\label{sect: Derivation of the Rényi Volume Forms}
As mentioned in Section \ref{sect: Rényi Volume Forms, Rényi-Priors}, the following technical derivation is based on the derivation of Takeuchi and Amari's $\alpha$-priors (see \autoref{eq: alpha parallel covolumes}) presented in Calin and Udriste \citep{Udriste}.

To begin with, consider the manifold constituted by the exponential family $\mathcal{P}_e$ (see \autoref{eq: exponential family}). 
The dual, parametrized Rényi-volume forms $\tensor[^{(\rho)}]{\omega}{}, \,\tensor[^{(\rho^*)}]{\omega}{}$ for $\mathcal{P}_e$ (which will be called $\tensor[^{(\rho)}]{\omega}{}_e$), must fulfill the following partial differential equations (PDEs) 
  \begin{equation}
    \tensor[^{(\rho)}]{\nabla}{}\tensor[^{(\rho)}]{\omega}{}_e  = 0\,, \quad \tensor[^{(\rho^*)}]{\nabla}{}   \tensor[^{(\rho^*)}]{\omega}{}_e  = 0\,,
  \end{equation}
Recall that on the manifold $(\mathcal{M}, \tensor[^\rho]{g}{})$ with the Rényi-metric $\tensor[^\rho]{g}{} = \rho \cdot \tensor[^{\text{F}}]{g}{}$, all top-dimensional forms are proportional, i.e. there exist two functions $f^\rho_e, f^{\rho^*}_e\in \mathcal{F}(\mathcal{M})$, s.t. 
  \begin{equation}
      \tensor[^{(\rho)}]{\omega}{}_e = f^\rho_e\, \dd V_\rho, \quad   \tensor[^{(\rho^*)}]{\omega}{}_e = f^{\rho^*}_e \, \dd V_\rho \,.\label{eq: definining eq for the e volume forms}
  \end{equation}
with $\dd V_\rho$ the Riemannian volume form with the metric tensor $\tensor[^\rho]{g}{}$. For now, we focus on $\tensor[^{(\rho^*)}]{\omega}{}_e$. Merely using the covariant derivative of a one-form, the rightmost expression in \autoref{eq: definining eq for the e volume forms} becomes 
\begin{align}
    \tensor[^{(\rho^*)}]{\nabla}{_{\partial_i}}  \tensor[^{(\rho^*)}]{\omega}{}_e &=(\partial_i f^{\rho^*}_e) \dd V_{\rho} + f^{\rho^*}_e \tensor[^{(\rho^*)}]{\nabla}{_{\partial_i}} \dd V_{\rho} \overset{!}{=}0 \\ 
    \leftrightarrow \; -\partial_i (\log f^{\rho^*}_e)\,\dd V_\rho &= \tensor[^{(\rho^*)}]{\nabla}{_{\partial_i}} \dd V_{\rho} \,.\label{eq: PDE for Rényi-parallel volume forms}
\end{align}
At this point, we plug \autoref{eq: rho* connection in terms of e and C} into \autoref{eq: definining eq for the e volume forms} to find 
\begin{align}
    -\partial_i (\log f^{\rho^*}_e) \dd V_{\rho} &\overset{\eqref{eq: rho* connection in terms of e and C}}{=} 2(1-\rho) \, \tensor[^{\text{LC}}]{\nabla}{} \,\dd V_\rho + 2\left(\rho-\frac{1}{2}\right) \tensor[^{(e)}]{\nabla}{}_{\partial_i}\; \dd V_\rho \\
    \leftrightarrow  \; -\partial_i (\log f^{\rho^*}_e) \dd V_{1} &= \underbrace{2(1-\rho) \, \tensor[^{\text{LC}}]{\nabla}{} \,\dd V_1}_{\sim \tensor[^{\text{LC}}]{\nabla}{}\,\dd V_1\, =\, 0} + 2\left(\rho-\frac{1}{2}\right) \tensor[^{(e)}]{\nabla}{}_{\partial_i}\; \dd V_1\label{eq: lhs rhs}
\end{align}
with $\dd V_1 = \sqrt{\det \tensor[^{\text{F}}]{g}{}} \;\dd x^1 \wedge \dots \wedge \dd x^n$. Evaluation of this expression in the tangent space's local induced basis leads to 
\begin{align}
-\partial_i (\log f^{\rho^*}_e) \dd V_{1}(\partial_1,\dots,\partial_n) &= -\partial_i (\log f^{\rho^*}_e) \, \sqrt{\det \tensor[^{\text{F}}]{g}{}} \,.
\label{eq: fe in terms of dV_1}
\end{align}
Furthermore, the covariant derivative of the alternating multilinear map $\dd V_1$ in \autoref{eq: lhs rhs} is carried out to find 
\begin{align}
    (\tensor[^{(e)}]{\nabla}{}_{\partial_i} \dd V_1)(\partial_1,\dots,\partial_n) 
    &=\partial_i \sqrt{\det \tensor[^{\text{F}}]{g}{}} - \sqrt{\det \tensor[^{\text{F}}]{g}{}}\; \tensor[^{(e)}]{\Gamma}{^{k}_{ki}} \,. 
\end{align} 
In summary, \autoref{eq: lhs rhs} thus gives 
\begin{align}
    -\partial_i (\log f^{\rho^*}_e) \, \sqrt{\det \tensor[^{\text{F}}]{g}{}} &= 2\left(\rho-\frac{1}{2}\right) \left(\partial_i \sqrt{\det \tensor[^{\text{F}}]{g}{}} - \sqrt{\det \tensor[^{\text{F}}]{g}{}}\; \tensor[^{(e)}]{\Gamma}{^{k}_{ki}}\right)\\
    \leftrightarrow \;  -\partial_i (\log f^{\rho^*}_e)  &= 2\left(\rho-\frac{1}{2}\right) \left(\partial_i \log \sqrt{\det \tensor[^{\text{F}}]{g}{}} - \tensor[^{(e)}]{\Gamma}{^{k}_{ki}}\right) \label{eq: f rho and det g and e connection}
\end{align} 
Combining these results with the flatness of the exponential family's statistical manifold with respect to $\tensor[^{(e)}]{\nabla}{}$, i.e. $\tensor[^{(e)}]{\Gamma}{} = 0$ on $\mathcal{P}_e$, one arrives at 
\begin{equation}
    \partial_i (\log f_e^{\rho^*}) = \partial_i \left(\log (\det \tensor[^{\text{F}}]{g}{})^{-\left(\rho-\frac{1}{2}\right)}\right) + 2 \left(\rho-\frac{1}{2}\right) \, \underbrace{\tensor[^{(e)}]{\Gamma}{^k_{ik}}}_{=0}\,,
\end{equation}
which ultimately gives 
\begin{align}
    0 &= \partial_i \log f_e^{\rho^*} ((\det \tensor[^{\text{F}}]{g}{})^{+ \left(\rho-\frac{1}{2}\right)}) \quad 
    \leftrightarrow \quad f_e^{\rho^*} = (\det \tensor[^{\text{F}}]{g}{})^{- \left(\rho-\frac{1}{2}\right)}\,,
\end{align}
and thus 
\begin{align}
      \tensor[^{(\rho^*)}]{\omega}{}_e &= (\det \tensor[^{\text{F}}]{g}{})^{- \left(\rho-\frac{1}{2}\right)} \; \dd V_\rho 
      &= \rho^{n/2}\cdot (\det \tensor[^{\text{F}}]{g}{})^{1-\rho} \;\dd x^1\wedge \dots \wedge \dd x^n\,.
\end{align}
By an analogous derivation, it is straightforward to see that  
\begin{equation}
    \tensor[^{(\rho)}]{\omega}{}_e = \rho^{n/2}\cdot (\det \tensor[^{\text{F}}]{g}{})^{\rho} \;\dd x^1\wedge \dots \wedge \dd x^n\,.
\end{equation}
The dual Rényi-parallel volume forms for the mixture family $\mathcal{P}_m$'s statistical manifold, $\tensor[^{(\rho)}]{\omega}{}_m\,,\omega^{\rho*}_m$ is are derived accordingly, see Section \ref{sect: the Results volume forms}.  
%%       \cite{<label>} ==> [1]

%%
% Example citation, See \cite{lamport94}.


\bibliographystyle{elsarticle-num} 
\bibliography{references}

\begin{thebibliography}{10}
\expandafter\ifx\csname url\endcsname\relax
  \def\url#1{\texttt{#1}}\fi
\expandafter\ifx\csname urlprefix\endcsname\relax\def\urlprefix{URL }\fi
\expandafter\ifx\csname href\endcsname\relax
  \def\href#1#2{#2} \def\path#1{#1}\fi

\bibitem{amari_information_2016}
S.-i. Amari, \href{http://link.springer.com/10.1007/978-4-431-55978-8}{Information Geometry and Its Applications}, Vol. 194 of Applied Mathematical Sciences, Springer Japan, 2016.
\newblock \href {https://doi.org/10.1007/978-4-431-55978-8} {\path{doi:10.1007/978-4-431-55978-8}}.
\newline\urlprefix\url{http://link.springer.com/10.1007/978-4-431-55978-8}

\bibitem{Shannonkhinchin0}
H.~Suyari, Generalization of shannon-khinchin axioms to nonextensive systems and the uniqueness theorem for the nonextensive entropy, IEEE Transactions on Information Theory 50~(8) (2004) 1783--1787.
\newblock \href {https://doi.org/10.1109/TIT.2004.831749} {\path{doi:10.1109/TIT.2004.831749}}.

\bibitem{Shannonkhinchin}
V.~M. Ilić, M.~S. Stanković, \href{https://www.sciencedirect.com/science/article/pii/S0378437114003732}{Generalized shannon–khinchin axioms and uniqueness theorem for pseudo-additive entropies}, Physica A: Statistical Mechanics and its Applications 411 (2014) 138--145.
\newblock \href {https://doi.org/10.1016/j.physa.2014.05.009} {\path{doi:10.1016/j.physa.2014.05.009}}.
\newline\urlprefix\url{https://www.sciencedirect.com/science/article/pii/S0378437114003732}

\bibitem{alphaparallelpriors}
J.~Takeuchi, S.~Amari, $\alpha$-parallel prior and its properties, IEEE Transactions on Information Theory 51~(3) (2005) 1011--1023.
\newblock \href {https://doi.org/10.1109/TIT.2004.842703} {\path{doi:10.1109/TIT.2004.842703}}.

\bibitem{Souza2016}
D.~C. de~Souza, R.~F. Vigelis, C.~C. Cavalcante, \href{https://www.mdpi.com/1099-4300/18/11/407}{Geometry induced by a generalization of rényi divergence}, Entropy 18~(11) (2016) 407.
\newblock \href {https://doi.org/10.3390/e18110407} {\path{doi:10.3390/e18110407}}.
\newline\urlprefix\url{https://www.mdpi.com/1099-4300/18/11/407}

\bibitem{van_erven_renyi_2014}
T.~Van~Erven, P.~Harremos, R{\'e}nyi divergence and kullback-leibler divergence, IEEE Transactions on Information Theory 60~(7) (2014) 3797--3820.

\bibitem{weylprior}
R.~Jiang, J.~Tavakoli, Y.~Zhao, Weyl prior and bayesian statistics, Entropy 22 (2020) 467.
\newblock \href {https://doi.org/10.3390/e22040467} {\path{doi:10.3390/e22040467}}.

\bibitem{Nielsen_2020}
F.~Nielsen, \href{https://doi.org/10.3390%2Fe22101100}{An elementary introduction to information geometry}, Entropy 22~(10) (2020) 1100.
\newblock \href {https://doi.org/10.3390/e22101100} {\path{doi:10.3390/e22101100}}.
\newline\urlprefix\url{https://doi.org/10.3390%2Fe22101100}

\bibitem{giesel_information_2021}
E.~Giesel, R.~Reischke, B.~M. Sch{\"a}fer, D.~Chia, \href{http://arxiv.org/abs/2005.01057}{Information geometry in cosmological inference problems}, JCAP 2021~(1) (2021) 005--005.
\newblock \href {http://arxiv.org/abs/2005.01057 [astro-ph, physics:gr-qc]} {\path{arXiv:2005.01057 [astro-ph, physics:gr-qc]}}, \href {https://doi.org/10.1088/1475-7516/2021/01/005} {\path{doi:10.1088/1475-7516/2021/01/005}}.
\newline\urlprefix\url{http://arxiv.org/abs/2005.01057}

\bibitem{rao1945information}
C.~R. Rao, Information and the accuracy attainable in the estimation of statistical parameters, Bulletin of the Calcutta Mathematical Society 37 (1945) 81--91.

\bibitem{original_fisher_information}
R.~A. Fisher, On the mathematical foundations of theoretical statistics, Philosophical transactions of the Royal Society of London. Series A, containing papers of a mathematical or physical character 222~(594-604) (1922) 309--368.

\bibitem{natural_gradient_descent}
M.~Rattray, D.~Saad, S.-i. Amari, Natural gradient descent for on-line learning, Physical review letters 81~(24) (1998) 5461.

\bibitem{rmhmc}
M.~Girolami, B.~Calderhead, \href{https://doi.org/10.1111/j.1467-9868.2010.00765.x}{Riemann manifold langevin and hamiltonian monte carlo methods}, Journal of the Royal Statistical Society Series B: Statistical Methodology 73~(2) (2011) 123--214.
\newblock \href {https://doi.org/10.1111/j.1467-9868.2010.00765.x} {\path{doi:10.1111/j.1467-9868.2010.00765.x}}.
\newline\urlprefix\url{https://doi.org/10.1111/j.1467-9868.2010.00765.x}

\bibitem{efron}
B.~Efron, \href{https://doi.org/10.1214/aos/1176343282}{{Defining the Curvature of a Statistical Problem (with Applications to Second Order Efficiency)}}, The Annals of Statistics 3~(6) (1975) 1189 -- 1242.
\newblock \href {https://doi.org/10.1214/aos/1176343282} {\path{doi:10.1214/aos/1176343282}}.
\newline\urlprefix\url{https://doi.org/10.1214/aos/1176343282}

\bibitem{jeffreysprior101}
H.~Jeffreys, \href{https://api.semanticscholar.org/CorpusID:19490929}{An invariant form for the prior probability in estimation problems}, Proceedings of the Royal Society of London. Series A. Mathematical and Physical Sciences 186 (1946) 453 -- 461.
\newline\urlprefix\url{https://api.semanticscholar.org/CorpusID:19490929}

\bibitem{manyfaces}
F.~Nielsen, The many faces of information geometry, Notices of the American Mathematical Society 69 (2022) 36--45.
\newblock \href {https://doi.org/10.1090/noti2403} {\path{doi:10.1090/noti2403}}.

\bibitem{eguchi1983divergenceinduced}
S.~Eguchi, Second order efficiency of minimum contrast estimators in a curved exponential family, The Annals of Statistics (1983) 793--803.

\bibitem{Bartelmann_2019}
M.~Bartelmann, \href{https://heiup.uni-heidelberg.de/catalog/book/534}{General Relativity}, Lecture Notes, Heidelberg University Publishing, 2019.
\newblock \href {https://doi.org/10.17885/heiup.534} {\path{doi:10.17885/heiup.534}}.
\newline\urlprefix\url{https://heiup.uni-heidelberg.de/catalog/book/534}

\bibitem{schaefer_tooltips_2022}
B.~M. Schäfer, \href{https://heiup.uni-heidelberg.de/catalog/book/1059}{Tooltips for Theoretical Physics: Concepts of Modern Theoretical Physics, Scales and Mathematical Tools}, Lecture Notes, Heidelberg University Publishing, 2022.
\newblock \href {https://doi.org/10.17885/heiup.1059} {\path{doi:10.17885/heiup.1059}}.
\newline\urlprefix\url{https://heiup.uni-heidelberg.de/catalog/book/1059}

\bibitem{kullback_leibler_original_definition}
S.~Kullback, R.~A. Leibler, On information and sufficiency, The annals of mathematical statistics 22~(1) (1951) 79--86.

\bibitem{chentsov1982statiscal}
N.~Chentsov, A.~M. Society, \href{https://books.google.de/books?id=kTJJtQAACAAJ}{Statistical Decision Rules and Optimal Inference}, Translations of mathematical monographs, American Mathematical Society, 1982.
\newline\urlprefix\url{https://books.google.de/books?id=kTJJtQAACAAJ}

\bibitem{amari_alpha_geometry_original_definition}
S.-i. Amari, Theory of information spaces: A differential geometrical foundation of statistics, Post RAAG Reports (1980).

\bibitem{renyi_original}
A.~R{\'e}nyi, et~al., On measures of information and entropy, in: Proceedings of the 4th Berkeley symposium on mathematics, statistics and probability, Vol.~1, 1961.

\bibitem{bhattacharyya_divergence_original_definition}
A.~Bhattacharyya, \href{http://www.jstor.org/stable/25047882}{On a measure of divergence between two multinomial populations}, Sankhy{\=a}: The Indian Journal of Statistics (1933-1960) 7~(4) (1946) 401--406.
\newline\urlprefix\url{http://www.jstor.org/stable/25047882}

\bibitem{Udriste}
O.~Calin, C.~Udriste, Geometric Modeling in Probability and Statistics, Springer Cham, 2014.
\newblock \href {https://doi.org/10.1007/978-3-319-07779-6} {\path{doi:10.1007/978-3-319-07779-6}}.

\bibitem{brown_admissible_1971}
L.~D. Brown, \href{http://projecteuclid.org/euclid.aoms/1177693318}{Admissible {Estimators}, {Recurrent} {Diffusions}, and {Insoluble} {Boundary} {Value} {Problems}}, The Annals of Mathematical Statistics 42~(3) (1971) 855--903.
\newblock \href {https://doi.org/10.1214/aoms/1177693318} {\path{doi:10.1214/aoms/1177693318}}.
\newline\urlprefix\url{http://projecteuclid.org/euclid.aoms/1177693318}

\bibitem{brown_heuristic_1979}
L.~D. Brown, \href{http://www.jstor.org/stable/2958667}{A {Heuristic} {Method} for {Determining} {Admissibility} of {Estimators}--{With} {Applications}}, The Annals of Statistics 7~(5) (1979) 960--994.
\newline\urlprefix\url{http://www.jstor.org/stable/2958667}

\bibitem{hartiganmaximumlikelihoodprior}
J.~A. Hartigan, \href{https://doi.org/10.1214/aos/1024691462}{{The maximum likelihood prior}}, The Annals of Statistics 26~(6) (1998) 2083 -- 2103.
\newblock \href {https://doi.org/10.1214/aos/1024691462} {\path{doi:10.1214/aos/1024691462}}.
\newline\urlprefix\url{https://doi.org/10.1214/aos/1024691462}

\bibitem{komaki_shrinkage_2001}
F.~Komaki, \href{https://doi.org/10.1093/biomet/88.3.859}{A shrinkage predictive distribution for multivariate {Normal} observables}, Biometrika 88~(3) (2001) 859--864.
\newblock \href {https://doi.org/10.1093/biomet/88.3.859} {\path{doi:10.1093/biomet/88.3.859}}.
\newline\urlprefix\url{https://doi.org/10.1093/biomet/88.3.859}

\bibitem{methodsamari}
S.-i. Amari, H.~Nagaoka, Methods of information geometry, Vol. 191, American Mathematical Soc., 2000.

\bibitem{Komaki_2006}
F.~Komaki, \href{http://dx.doi.org/10.1214/009053606000000010}{Shrinkage priors for bayesian prediction}, The Annals of Statistics 34~(2) (Apr. 2006).
\newblock \href {https://doi.org/10.1214/009053606000000010} {\path{doi:10.1214/009053606000000010}}.
\newline\urlprefix\url{http://dx.doi.org/10.1214/009053606000000010}

\bibitem{komaki2015}
F.~Komaki, \href{https://doi.org/10.1214/14-BA886}{{Asymptotic Properties of Bayesian Predictive Densities When the Distributions of Data and Target Variables are Different}}, Bayesian Analysis 10~(1) (2015) 31 -- 51.
\newblock \href {https://doi.org/10.1214/14-BA886} {\path{doi:10.1214/14-BA886}}.
\newline\urlprefix\url{https://doi.org/10.1214/14-BA886}

\bibitem{Kass1996TheSO}
R.~E. Kass, L.~A. Wasserman, \href{https://api.semanticscholar.org/CorpusID:53645083}{The selection of prior distributions by formal rules}, Journal of the American Statistical Association 91 (1996) 1343--1370.
\newline\urlprefix\url{https://api.semanticscholar.org/CorpusID:53645083}

\bibitem{Matsuzoe2015InformationGO}
H.~Matsuzoe, \href{https://doi.org/10.1063/1.4905989}{Information geometry of bayesian statistics}, AIP Conference Proceedings 1641~(1) (2015) 279--286.
\newblock \href {https://doi.org/10.1063/1.4905989} {\path{doi:10.1063/1.4905989}}.
\newline\urlprefix\url{https://doi.org/10.1063/1.4905989}

\bibitem{nakahara2003geometry}
M.~Nakahara, \href{https://books.google.de/books?id=cH-XQB0Ex5wC}{Geometry, Topology and Physics, Second Edition}, Graduate student series in physics, Taylor \& Francis, 2003.
\newline\urlprefix\url{https://books.google.de/books?id=cH-XQB0Ex5wC}

\bibitem{lee2019introduction}
J.~Lee, \href{https://books.google.de/books?id=UIPltQEACAAJ}{Introduction to Riemannian Manifolds}, Graduate Texts in Mathematics, Springer International Publishing, 2019.
\newline\urlprefix\url{https://books.google.de/books?id=UIPltQEACAAJ}

\bibitem{hartiganinvariantpriordistributions}
J.~Hartigan, \href{http://www.jstor.org/stable/2238537}{Invariant prior distributions}, The Annals of Mathematical Statistics 35~(2) (1964) 836--845.
\newline\urlprefix\url{http://www.jstor.org/stable/2238537}

\bibitem{hartigan1965}
J.~A. Hartigan, \href{https://doi.org/10.1214/aoms/1177699988}{{The Asymptotically Unbiased Prior Distribution}}, The Annals of Mathematical Statistics 36~(4) (1965) 1137 -- 1152.
\newblock \href {https://doi.org/10.1214/aoms/1177699988} {\path{doi:10.1214/aoms/1177699988}}.
\newline\urlprefix\url{https://doi.org/10.1214/aoms/1177699988}

\end{thebibliography}




% %% For numbered reference style

%   \textit{\LaTeX: a document preparation system},
%   Addison Wesley, Massachusetts,

% \end{thebibliography}
\end{document}